\documentclass[11pt]{article}
\usepackage{amsmath}
\usepackage{amssymb}
\usepackage{amscd}
\usepackage{xspace}
\usepackage{verbatim}
\newcommand{\mysection}[1]{
\section{#1}\setcounter{equation}{0}}
\title{\bf Local and global properties of solutions of  heat equation with superlinear absorption}
\author{{\bf Tai Nguyen Phuoc}\quad
 {\bf Laurent V\'eron}\\[2mm]
{\small Laboratoire de Math\'ematiques et Physique Th\'eorique, }\\
{\small  Universit\'e Fran\c{c}ois Rabelais,  Tours,  FRANCE}}

\date{}
\begin{document}
 \maketitle


\newcommand{\txt}[1]{\;\text{ #1 }\;}
\newcommand{\tbf}{\textbf}
\newcommand{\tit}{\textit}
\newcommand{\tsc}{\textsc}
\newcommand{\trm}{\textrm}
\newcommand{\mbf}{\mathbf}
\newcommand{\mrm}{\mathrm}
\newcommand{\bsym}{\boldsymbol}
\newcommand{\scs}{\scriptstyle}
\newcommand{\sss}{\scriptscriptstyle}
\newcommand{\txts}{\textstyle}
\newcommand{\dsps}{\displaystyle}
\newcommand{\fnz}{\footnotesize}
\newcommand{\scz}{\scriptsize}
\newcommand{\be}{\begin{equation}}
\newcommand{\bel}[1]{\begin{equation}\label{#1}}
\newcommand{\ee}{\end{equation}}
\newcommand{\eqnl}[2]{\begin{equation}\label{#1}{#2}\end{equation}}
\newcommand{\barr}{\begin{eqnarray}}
\newcommand{\earr}{\end{eqnarray}}
\newcommand{\bars}{\begin{eqnarray*}}
\newcommand{\ears}{\end{eqnarray*}}
\newcommand{\nnu}{\nonumber \\}
\newtheorem{subn}{\name}
\renewcommand{\thesubn}{}
\newcommand{\bsn}[1]{\def\name{#1}\begin{subn}}
\newcommand{\esn}{\end{subn}}
\newtheorem{sub}{\name}[section]
\newcommand{\dn}[1]{\def\name{#1}}   
\newcommand{\bs}{\begin{sub}}
\newcommand{\es}{\end{sub}}
\newcommand{\bsl}[1]{\begin{sub}\label{#1}}
\newcommand{\bth}[1]{\def\name{Theorem}
\begin{sub}\label{t:#1}}
\newcommand{\blemma}[1]{\def\name{Lemma}
\begin{sub}\label{l:#1}}
\newcommand{\bcor}[1]{\def\name{Corollary}
\begin{sub}\label{c:#1}}
\newcommand{\bdef}[1]{\def\name{Definition}
\begin{sub}\label{d:#1}}
\newcommand{\bprop}[1]{\def\name{Proposition}
\begin{sub}\label{p:#1}}
\newcommand{\R}{\eqref}
\newcommand{\rth}[1]{Theorem~\ref{t:#1}}
\newcommand{\rlemma}[1]{Lemma~\ref{l:#1}}
\newcommand{\rcor}[1]{Corollary~\ref{c:#1}}
\newcommand{\rdef}[1]{Definition~\ref{d:#1}}
\newcommand{\rprop}[1]{Proposition~\ref{p:#1}}
\newcommand{\BA}{\begin{array}}
\newcommand{\EA}{\end{array}}
\newcommand{\BAN}{\renewcommand{\arraystretch}{1.2}
\setlength{\arraycolsep}{2pt}\begin{array}}
\newcommand{\BAV}[2]{\renewcommand{\arraystretch}{#1}
\setlength{\arraycolsep}{#2}\begin{array}}
\newcommand{\BSA}{\begin{subarray}}
\newcommand{\ESA}{\end{subarray}}
\newcommand{\BAL}{\begin{aligned}}
\newcommand{\EAL}{\end{aligned}}
\newcommand{\BALG}{\begin{alignat}}
\newcommand{\EALG}{\end{alignat}}
\newcommand{\BALGN}{\begin{alignat*}}
\newcommand{\EALGN}{\end{alignat*}}
\newcommand{\note}[1]{\textit{#1.}\hspace{2mm}}
\newcommand{\Proof}{\note{Proof}}
\newcommand{\qeda}{\hspace{10mm}\hfill $\square$}
\newcommand{\qed}{\\
${}$ \hfill $\square$}
\newcommand{\Remark}{\note{Remark}}
\newcommand{\modin}{$\,$\\[-4mm] \indent}
\newcommand{\forevery}{\quad \forall}
\newcommand{\set}[1]{\{#1\}}
\newcommand{\setdef}[2]{\{\,#1:\,#2\,\}}
\newcommand{\setm}[2]{\{\,#1\mid #2\,\}}
\newcommand{\mt}{\mapsto}
\newcommand{\lra}{\longrightarrow}
\newcommand{\lla}{\longleftarrow}
\newcommand{\llra}{\longleftrightarrow}
\newcommand{\Lra}{\Longrightarrow}
\newcommand{\Lla}{\Longleftarrow}
\newcommand{\Llra}{\Longleftrightarrow}
\newcommand{\warrow}{\rightharpoonup}
\newcommand{
\paran}[1]{\left (#1 \right )}
\newcommand{\sqbr}[1]{\left [#1 \right ]}
\newcommand{\curlybr}[1]{\left \{#1 \right \}}
\newcommand{\abs}[1]{\left |#1\right |}
\newcommand{\norm}[1]{\left \|#1\right \|}
\newcommand{
\paranb}[1]{\big (#1 \big )}
\newcommand{\lsqbrb}[1]{\big [#1 \big ]}
\newcommand{\lcurlybrb}[1]{\big \{#1 \big \}}
\newcommand{\absb}[1]{\big |#1\big |}
\newcommand{\normb}[1]{\big \|#1\big \|}
\newcommand{
\paranB}[1]{\Big (#1 \Big )}
\newcommand{\absB}[1]{\Big |#1\Big |}
\newcommand{\normB}[1]{\Big \|#1\Big \|}

\newcommand{\thkl}{\rule[-.5mm]{.3mm}{3mm}}
\newcommand{\thknorm}[1]{\thkl #1 \thkl\,}
\newcommand{\trinorm}[1]{|\!|\!| #1 |\!|\!|\,}
\newcommand{\bang}[1]{\langle #1 \rangle}
\def\angb<#1>{\langle #1 \rangle}
\newcommand{\vstrut}[1]{\rule{0mm}{#1}}
\newcommand{\rec}[1]{\frac{1}{#1}}
\newcommand{\opname}[1]{\mbox{\rm #1}\,}
\newcommand{\supp}{\opname{supp}}
\newcommand{\dist}{\opname{dist}}
\newcommand{\myfrac}[2]{{\displaystyle \frac{#1}{#2} }}
\newcommand{\myint}[2]{{\displaystyle \int_{#1}^{#2}}}
\newcommand{\mysum}[2]{{\displaystyle \sum_{#1}^{#2}}}
\newcommand {\dint}{{\displaystyle \int\!\!\int}}
\newcommand{\q}{\quad}
\newcommand{\qq}{\qquad}
\newcommand{\hsp}[1]{\hspace{#1mm}}
\newcommand{\vsp}[1]{\vspace{#1mm}}
\newcommand{\ity}{\infty}
\newcommand{\prt}{\partial}
\newcommand{\sms}{\setminus}
\newcommand{\ems}{\emptyset}
\newcommand{\ti}{\times}
\newcommand{\pr}{^\prime}
\newcommand{\ppr}{^{\prime\prime}}
\newcommand{\tl}{\tilde}
\newcommand{\sbs}{\subset}
\newcommand{\sbeq}{\subseteq}
\newcommand{\nind}{\noindent}
\newcommand{\ind}{\indent}
\newcommand{\ovl}{\overline}
\newcommand{\unl}{\underline}
\newcommand{\nin}{\not\in}
\newcommand{\pfrac}[2]{\genfrac{(}{)}{}{}{#1}{#2}}

\def\ga{\alpha}     \def\gb{\beta}       \def\gg{\gamma}
\def\gc{\chi}       \def\gd{\delta}      \def\ge{\epsilon}
\def\gth{\theta}                         \def\vge{\varepsilon}
\def\gf{\phi}       \def\vgf{\varphi}    \def\gh{\eta}
\def\gi{\iota}      \def\gk{\kappa}      \def\gl{\lambda}
\def\gm{\mu}        \def\gn{\nu}         \def\gp{\pi}
\def\vgp{\varpi}    \def\gr{\rho}        \def\vgr{\varrho}
\def\gs{\sigma}     \def\vgs{\varsigma}  \def\gt{\tau}
\def\gu{\upsilon}   \def\gv{\vartheta}   \def\gw{\omega}
\def\gx{\xi}        \def\gy{\psi}        \def\gz{\zeta}
\def\Gg{\Gamma}     \def\Gd{\Delta}      \def\Gf{\Phi}
\def\Gth{\Theta}
\def\Gl{\Lambda}    \def\Gs{\Sigma}      \def\Gp{\Pi}
\def\Gw{\Omega}     \def\Gx{\Xi}         \def\Gy{\Psi}

\def\CS{{\mathcal S}}   \def\CM{{\mathcal M}}   \def\CN{{\mathcal N}}
\def\CR{{\mathcal R}}   \def\CO{{\mathcal O}}   \def\CP{{\mathcal P}}
\def\CA{{\mathcal A}}   \def\CB{{\mathcal B}}   \def\CC{{\mathcal C}}
\def\CD{{\mathcal D}}   \def\CE{{\mathcal E}}   \def\CF{{\mathcal F}}
\def\CG{{\mathcal G}}   \def\CH{{\mathcal H}}   \def\CI{{\mathcal I}}
\def\CJ{{\mathcal J}}   \def\CK{{\mathcal K}}   \def\CL{{\mathcal L}}
\def\CT{{\mathcal T}}   \def\CU{{\mathcal U}}   \def\CV{{\mathcal V}}
\def\CZ{{\mathcal Z}}   \def\CX{{\mathcal X}}   \def\CY{{\mathcal Y}}
\def\CW{{\mathcal W}} \def\CQ{{\mathcal Q}}
\def\BBA {\mathbb A}   \def\BBb {\mathbb B}    \def\BBC {\mathbb C}
\def\BBD {\mathbb D}   \def\BBE {\mathbb E}    \def\BBF {\mathbb F}
\def\BBG {\mathbb G}   \def\BBH {\mathbb H}    \def\BBI {\mathbb I}
\def\BBJ {\mathbb J}   \def\BBK {\mathbb K}    \def\BBL {\mathbb L}
\def\BBM {\mathbb M}   \def\BBN {\mathbb N}    \def\BBO {\mathbb O}
\def\BBP {\mathbb P}   \def\BBR {\mathbb R}    \def\BBS {\mathbb S}
\def\BBT {\mathbb T}   \def\BBU {\mathbb U}    \def\BBV {\mathbb V}
\def\BBW {\mathbb W}   \def\BBX {\mathbb X}    \def\BBY {\mathbb Y}
\def\BBZ {\mathbb Z}

\def\GTA {\mathfrak A}   \def\GTB {\mathfrak B}    \def\GTC {\mathfrak C}
\def\GTD {\mathfrak D}   \def\GTE {\mathfrak E}    \def\GTF {\mathfrak F}
\def\GTG {\mathfrak G}   \def\GTH {\mathfrak H}    \def\GTI {\mathfrak I}
\def\GTJ {\mathfrak J}   \def\GTK {\mathfrak K}    \def\GTL {\mathfrak L}
\def\GTM {\mathfrak M}   \def\GTN {\mathfrak N}    \def\GTO {\mathfrak O}
\def\GTP {\mathfrak P}   \def\GTR {\mathfrak R}    \def\GTS {\mathfrak S}
\def\GTT {\mathfrak T}   \def\GTU {\mathfrak U}    \def\GTV {\mathfrak V}
\def\GTW {\mathfrak W}   \def\GTX {\mathfrak X}    \def\GTY {\mathfrak Y}
\def\GTZ {\mathfrak Z}   \def\GTQ {\mathfrak Q}

\font\Sym= msam10 
\def\SYM#1{\hbox{\Sym #1}}
\newcommand{\bdw}{\prt\Gw\xspace}
\medskip
\begin{abstract}
We study the limit, when $k\to\infty$ of the
solutions of $ \prt_tu-\Delta u+f(u)=0$ in $\BBR^N\ti(0,\infty)$ with initial data $k\gd$, when $f$ is a positive superlinear increasing function. We prove that there exist essentially three types of possible behaviour according $f^{-1}$ and $F^{-1/2}$ belong or not to $L^1(1,\infty)$, where $F(t)=\int_0^t f(s)ds$. We use these results for providing a new and more general construction of the initial trace and some uniqueness and non-uniqueness results for solutions with unbounded initial data.
\end{abstract}

\tableofcontents
\mysection{Introduction}

In this article we investigate some local and global properties of solutions of a class of semilinear heat equations
\bel{A1} \prt_t u -\Gd u + f(u) = 0 \ee
in $Q_\ity:={\BBR}^N \ti (0,\ity)$ $(N \geq 2)$ where $f:\BBR_+\mapsto\BBR_+$ is continuous, nondecreasing and positive on $(0,\infty)$, vanishes at $0$ and tends to infinity at infinity.  As a {\it model equation} we shall consider the following nonlinear term, with $\ga>0$,
\bel{A2} 
\prt_t u -\Gd u + u\ln^\ga(u+1) = 0, \ee
which points out all the delicate features of {\it weakly superlinear absorption}. By opposition, for power-like absorption $f(u)=|u|^{\gb}u$ with $\gb>0$ much is known about the structure of the set of solutions. The local and asymptotic behaviour of solutions is strongly linked to the existence of a self-similar solutions under the form 
\bel{A3} 
u(x,t)= t^{-1/\gb}w(x\sqrt t). 
\ee
In this case the critical exponent $\gb_c=2/N$ plays a fundamental role in the  description of isolated singularities and the study of the initial trace. This is due to the fact that, for $0<\gb<\gb_c$, there exists a positive self-similar solution with an isolated singularity at $(0,0)$ and vanishing on $\BBR^N\setminus\{0\}\ti\{0\}$, while no such solution exists when $\gb\geq \gb_c$ and more generally, no solution with isolated singularities.

In the case of $(\ref{A2})$, no self-similar structure exists. There is no critical exponent corresponding to isolated singularities since there always exist such singular solutions. Actually, for any $k>0$ there exists a unique 
$u=u_k\in C(\overline {Q_\infty}\setminus \{(0,0)\})\cap C^{2,1}(Q_\infty)$  solution of 
\bel{A4} 
\left\{\BA {ll}
\prt_t u-\Gd u+u\ln^\ga(u+1)=0\qquad&\text {in }Q_\infty\\[2mm]\phantom{-------....}
u(.,x)=k\gd_0\qquad&\text {in }\CD'(\BBR^N).
\EA\right.\ee
There are two critical values for $\ga$: $\ga=1$ and $\ga=2$, the explanation of which comes from the study of the two singular problems
\bel{A5} 
\left\{\BA {ll}
\gf'+\gf\ln^\ga(\gf+1)&=0\qquad\text {in }(0,\infty)\\[2mm]\phantom{------}
\gf(0)&=\infty,
\EA\right.\ee
and, for any $\ge>0$, 
\bel{A6} \left\{\BA {ll}
-\Gd\psi+\psi\ln^\ga(\psi+1)&=0\qquad\text {in }\BBR^N\setminus B_\ge\\[2mm]\phantom{;;;;----}
{\displaystyle \lim_{|x|\to \ge}}\psi(x)&=\infty,
\EA\right.\ee
where $B_\ge:=\{x\in\BBR^N:|x|<\ge\}$. When it exists, the solution $\gf_\infty$ of $(\ref{A5})$ is given implicitely by
\bel{A7} 
\myint{\gf_\infty(t)}{\infty}\myfrac{ds}{s\ln^\ga(s+1)}=t\qquad\forall t>0,
\ee
and such a formula is valid if and only if $\ga>1$. For problem $(\ref{A6})$ an explicit expression of the solution  is not valid, but this solution exists if and only if  $\ga>2$; in this case of the Keller-Osserman condition (see $(\ref{A12})$ below) holds. \smallskip

Having in mind this model we study $(\ref{A1})$ assuming the {\it weak singularity condition} on $f$:
\bel{A8} 
\myint{1}{\infty}s^{-2-\frac{2}{N}}f(s)ds<\infty.
\ee
\bprop{exist-fundamental}
Assume $(\ref{A8})$ holds. Then for any $k>0$, there exists a unique solution $u:=u_k$ to 
\bel{A9} \left\{ \BA{lll} \prt_t u -\Gd u +f(u) &= 0 \qq &\text{in } Q_\ity \\ \phantom{-----}
u(.,0) &= k\gd_0 &\text{in } \CD'(\BBR^N).
\EA \right. \ee 
Furthermore, if $\psi_n$ is a sequence of positive integrable functions converging to $k\gd$ in the weak-star topology, then the sequence $u_{\psi_n}$ of solutions of $(\ref{A1})$ in $Q_\infty$ with initial data $\psi_n$ converges to $u_{k\gd}$, locally uniformly.
\es 
Another important condition on $f$ is
\bel{A10} 
\myint{1}{\infty}\myfrac{ds}{f(s)}<\infty.
\ee
Under assumption $(\ref{A10})$ there exists a solution $\gf:=\gf_\ity$ to
\bel{A11} 
\left\{\BA {ll}
\gf'+f(\gf)&=0\qquad\text {in }(0,\infty)\\[2mm]\phantom{--;;}
\gf(0)&=\infty.
\EA\right.\ee
The function $\gf_\infty$ is the maximal solution of (\ref{A11}) and it it explicited by a formula similar to (\ref{A7}) in which $s\ln^\ga(s+1)$ is replaced by $f(s)$. 

The next important condition on $f$ we shall encounter is the Keller-Osserman condition, i.e.
\bel{A12} 
\myint{1}{\infty}\myfrac{ds}{\sqrt{F(s)}}<\infty,
\ee
where 
\bel{A13} F(s)=\myint{0}{s}f(\gs)d\gs, \q \forall s \in [1,\ity).
\ee
If $(\ref{A12})$ is satisfied, by \cite[Theorem III]{Ke} for any $\ge>0$ there exists a maximal solution $\psi:=\psi_\ge$ to
\bel{A14} \left\{\BA {ll}
-\Gd\psi+ f(\psi)&=0\qquad\text {in }\BBR^N\setminus B_\ge\\[2mm]\phantom{--}
{\displaystyle \lim_{|x|\to \ge}}\psi(x)&=\infty.
\EA\right.\ee

Assumptions (\ref{A10}) and (\ref{A13}) which are simultaneously satisfied in the case of a power like absorption, but not in our model case, are the Ariane shred which illuminates the structure of the set of solutions of (\ref{A1}), in particular in view of the initial trace problem.

The first question we consider is the study of  the limit of $u_k$ when $k\to\infty$. 
This question is natural since $k\mapsto u_k$ is increasing. In order to treat it, we need some additional conditions. \smallskip
 
\noindent (C1)- The function $s \mapsto \myfrac{f(s)}{s}$ is increasing on $(0,\ity)$ and satisfies 
$${\displaystyle \lim_{s \to 0}}\myfrac{f(s)}{s}=0 \q \text{and} \q {\displaystyle \lim_{s \to \ity}}\myfrac{f(s)}{s}=\ity.$$
  
\noindent (C2)- The function $f$ is convex on $(0,\ity)$.\smallskip
 
\noindent (C3)-  If ${\displaystyle \liminf_{s \to \infty}}\,f(s)/(s\ln^\ga s)=0,\forall{\ga >2}$, then there exists $\gb \in (1,2]$ such that 
\[ \limsup_{s \to \infty}\myfrac{f(s)}{s\ln^\gb s}<\ity. \] 
In the second section, we prove the following results.\medskip
\bth{ALA} Assume the conditions $(C1)$ and $(C3)$ hold. If $f$ satisfies
\bel{A15} 
\myint{1}{\infty}\myfrac{ds}{f(s)}=\infty,
\ee
then the solutions $u_k$ of $(\ref{A9})$ satisfy ${\displaystyle \lim_{k\to\infty}}u_k(x,t)=\infty$ for every $(x,t) \in Q_\ity$.
\es 
\bth{TTA} Assume the conditions $(C1)-(C3)$ hold. If $f$ satisfies $(\ref{A10})$ and
\bel{A16} 
\myint{1}{\infty}\myfrac{ds}{\sqrt{F(s)}}=\infty
\ee
where $F$ is defined in $(\ref{A13})$, then the solutions $u_k$ of $(\ref{A9})$ satisfy ${\displaystyle \lim_{k\to\infty}}u_k(x,t)=\gf_\infty(t)$ for every $(x,t) \in Q_\ity$, where $\gf_\ity$ is the solution of $(\ref{A11})$.
\es
We denote by $\CU_0$ the set of positive solutions $u$ of $(\ref{A1})$ in $Q_\infty$, which are continuous in $\overline {Q_\infty}\setminus \{(0,0)\}$, vanish on the set $\{(x,0):x\neq 0\}$ and satisfies
\bel{A17} 
\lim_{t\to 0}\myint{B_\ge}{}u(x,t)dx=\infty
\ee
for any $\ge>0$.
\bth{min-solution} Assume $f$ satisfies $(\ref{A8})$, $(\ref{A12})$ and (C2). Then ${\unl U:=\displaystyle \lim_{k\to\infty}}u_k$ is the minimal element of $\CU_0$.
\es

In the third section we study the set of positive and locally bounded solutions of $(\ref{A1})$ in $Q_\infty$. This set differs considerably according the assumption on $f$. This is due to the properties of the radial solutions of the associated stationnary equation
\bel{A18} 
-\Gd w+f(w)=0\qquad\text{in }\BBR^N.
\ee
The next result is based upon the Picard-Lipschitz fixed point theorem and a result of Vazquez and V\'eron \cite{VaVe}.
\bprop{Ex1} Assume $(\ref{A16})$ holds. For any $a>0$, there exists a unique positive function $w:=w_a \in C^2([0,\ity))$ to the problem
\bel{A19}  \left\{ \BA{ll} 
-w \ppr-\myfrac{N-1}{r}w\pr+f(w)&=0 \qq \text{in } \BBR_+ \\[1mm]\phantom{-w \ppr-\myfrac{N-1}{r}w\pr+,}
w\pr(0)&=0 \\\phantom{-w \ppr-\myfrac{N-1}{r}w\pr,,,}
w(0)&=a.
\EA \right. \ee
\es
A striking consequence of the existence of such solutions is the following non-uniqueness result. 

\bth{non-unique} Assume $f$ satisfies $(\ref{A10})$ and $(\ref{A16})$. Then for any $u_0\in C(\BBR^N)$ satisfying, for some $b>a>0$, $w_a(x)\leq u_0(x)\leq w_b(x)$ $\forall x \in \BBR^N$, there exist two solutions $\underline u, \overline u\in C(\overline {Q_\infty})$ of $(\ref{A1})$ with initial value $u_0$. They satisfy respectively
\bel{A20} 
0\leq \underline u(x,t)\leq \min\{w_b(x),\phi_\infty(t)\}\qquad\forall (x,t)\in Q_\infty,
\ee
thus ${\displaystyle \lim_{t\to\infty}}\underline u(x,t)=0$, uniformly with respect to $x \in \BBR^N$, and
\bel{A21} 
w_a(x)\leq \overline u(x,t)\leq w_b(x)\qquad\forall (x,t)\in Q_\infty,
\ee
thus ${\displaystyle \lim_{|x|\to\infty}}\overline u(x,t)=\infty$, uniformly with respect to $t\geq 0$. \es 

The next theorem shows that if two solutions of $(\ref{A1})$ have the same initial data and the same asymptotic behaviour as $\abs x \to \ity$ then they coincide.  
\bth{unique} Assume $f$ satisfies $(C1)$ and $(\ref{A16})$. Let $u, \tl u \in C(\ovl Q_\ity) \cap C^{2,1}(Q_\ity)$ be two positive solutions of $(\ref{A1})$ with initial data $u_0$. If for any $\ge >0$,
\bel{A22} u(x,t)-\tl u(x,t) = o(w_\ge(\abs x)) \text{ as } x \to \ity \ee
locally uniformly with respect to $t\geq 0$, then $u=\tl u$.
\es 
On the contrary, if the Keller-Osserman condition holds, a continuous solution is uniquely determined by the positive initial value $u_0 \in C(\BBR^N)$, and uniqueness still holds if $C(\BBR^N)$ is replaced by $\mathfrak M_+(\BBR^N)$.
\bth{unique-KO} Assume $f$ satisfies $(\ref{A12})$ and $(C2)$. Then\smallskip

\noindent (i)  For any nonnegative function $u_0 \in C(\BBR^N)$ there exists a unique nonnegative solution $u \in C(\ovl Q_\ity)$ of $(\ref{A1})$ in $Q_\ity$ with initial value $u_0$. \smallskip

\noindent (ii)  For any for any nonnegative measure $\gm \in \mathfrak M(\BBR^N)$, there exists at most one nonnegative solution $u \in C(Q_\ity)$ of $(\ref{A1})$ in $Q_\ity$ such that $f(u)\in L^1_{loc}(\ovl Q_\ity)$ satisfying
\bel{A22'} \lim_{t\to 0}\myint{\BBR^N}{} u(x,t)\gz (x)dx=\myint{\BBR^N}{} \gz (x)d\gm(x)\quad \forall\gz\in C_c(\BBR^N).
\ee
\es \medskip  

In the last section we use the tools studied in the previous sections to develop a {\it new construction of the initial trace} of locally bounded positive solutions of $(\ref{A1})$ in $Q_\infty$. By opposition to the power-like case \cite {MV1}, where the initial trace was constructed by duality arguments based upon H\"older inequality and delicate choice of test functions, our new method has the advantage of being based only on maximum principle, using either the Keller-Osserman condition, if (\ref{A16}) holds, or the asymptotics of the $u_k$ if (\ref{A16}) does not hold. We first prove

\bprop{reg} Let $u\in C^{2,1}(Q_\infty)$ be a positive solution of $(\ref{A1})$ in $Q_\infty$. The set $\CR(u)$ of the points $z\in \BBR^N$ such that there exists an open ball $B_r(z)$ such that $u, f(u)\in L^1(Q_{T}^{B_r(z)})$ is an open subset. Furthermore there exists a positive Radon measure $\gm:=\gm(u)$ on $\CR(u)$ such that
\bel{A23} \lim_{t \to 0}\myint{\CR(u)}{} u(x,t)\gz(x)dx=\myint{\CR(u)}{} \gz(x)d\gm(x)\qq\forall \gz\in C_c(\CR(u)).\ee
\es
Due to \rprop{reg}, we introduce the definition of the initial trace.
\bdef{Reg} The couple $(\CS(u),\gm)$ where $\CS(u)=\BBR^N\sms \CR(u)$ is called the initial trace of $u$ in $\Gw$ and will be denoted by $tr_{\BBR^N}(u)$. The set $\CR(u)$ is called the regular set of the initial trace of $u$ and the measure $\gm$ the regular part of the initial trace. The set $\CS(u)$ is closed and is called the singular part of the initial trace of $u$.  
\es
The initial trace can also be represented by a positive, outer regular Borel measure, not necessary locally bounded. The space of these measures on $\BBR^N$ will be denoted by $\CB_+^\text{reg}(\BBR^N)$. If for every open subset $A \sbs \BBR^N$ we denote by $\GTM_+(A)$ the space of positive Radon measures on $A$, there is a one-to-one correspondence between $\CB_+^\text{reg}(\BBR^N)$ and the set of couples:
\bel{A24} CM_+(\BBR^N)=\left\{(\CS,\gm): \CS \sbs \BBR^N \text{  closed}, \gm \in \GTM_+(\CR) \text{ with } \CR=\BBR^N \sms \CS\right\}. \ee
The Borel measure $\gn \in \CB_+^\text{reg}(\BBR^N)$ corresponding to a couples $(\CS,\gm) \in CM_+(\BBR^N)$ is given by
\bel{A25} \gn(A)= \left\{ \BA{ll} \ity \qq &\text{if } A \cap \CS \ne \ems  \\
\gm(A) &\text{if } A \subseteq \CS, \EA \right. \forevery A \sbs \BBR^N, A \text{ Borel. } \ee
If $u$ is a solution of $(\ref{A1})$, we shall use the notation $tr_{\BBR^N}(u)$ (resp. $Tr_{\BBR^N}(u)$) for the trace considered as an element of $CM_+(\BBR^N)$ (resp. $\CB_+^\text{reg}(\BBR^N)$). \medskip

We consider the case when the Keller-Osserman holds.

\bth{tr+KO} Assume $f$ is nondecreasing and satisfies $(\ref{A12})$. If $u\in C^{2,1}(Q_\infty)$ is a positive solution of $(\ref{A1})$, it possesses an initial trace $\gn\in\CB_+^{reg}(\BBR^N)$.
\es

Furthermore, the following theorem deals with the existence of the maximal solution and the minimal solution of $(\ref{A1})$ with a given initial trace $(\CS,\gm) \in CM_+(\BBR^N)$.
\bth{KO-existence}  Assume  $f$ is nondecreasing and satisfies $(\ref{A12})$, $(\ref{A8})$ and $(C2)$. Then for any $(\CS,\gm) \in CM_+(\BBR^N)$ there exist a maximal solution $\ovl u_{\CS,\gm}$ and a minimal solution $\unl u_{\CS,\gm}$ of $(\ref{A1})$ in $Q_\infty$, with initial trace $(\CS,\gm)$, in the following sense:
\bel{A26} \unl u_{\CS,\gm} \leq v \leq \ovl u_{\CS,\gm} \ee
for every positive solution $v \in C^{2,1}(Q_{\infty})$ of $(\ref{A1})$ in $Q_\infty$ such that $tr_{\BBR^N}(v)=(\CS,\gm)$. 
\es \medskip

If the Keller-Osserman does not holds, we obtain the following results which depend upon $\displaystyle\lim_{k\to\infty}u_k$ is equal to $\phi_\infty$ or is infinite (we recall that $u_k$ is the solution of $(\ref{A9})$).

\bth{tr+J-fin} Assume $(\ref{A8})$, $(\ref{A10})$ and $(\ref{A16})$ are verified and $\displaystyle\lim_{k\to\infty}u_k=\gf_\infty$. If $u$ is a positive solution of $(\ref{A1})$ in $Q_{\infty}$, it possesses an initial trace which is either the Borel measure $\gn_\infty$ which satisfies $\gn_\infty(\CO)=\infty$ for any non-empty open subset $\CO\subset \BBR^N$, or is a positive Radon measure $\gm$ on $\BBR^N$.
This result holds in particular if $(C1)$ and $(C3)$ hold.
\es

A consequence of \rth {tr+J-fin} which is worth mentioning is the following. 
\bprop {exit-2} Under the assumptions of \rth {tr+J-fin}, for any $b>0$ there exists a positive solution $u\in C(Q_\infty)$ of $(\ref{A1})$ in $(\ref{A1})$ satisfying
  \bel{J-10} 
  \max\{\phi_\infty(t);w_b(|x|)\}\leq  u(x,t)\leq \phi_\infty(t)+w_b(|x|)\qquad \forall (x,t)\in Q_\infty.
  \ee
Consequently there exist infinitely many positive solutions of  $(\ref{A1})$ with initial trace $\gn_\infty$. Furthermore
$\phi_\infty$ is the smallest of all these solutions. 
\es

\bth{tr+J-infin} Assume $f$ satisfies $(\ref{A8})$, $(\ref{A15})$, $(\ref{A16})$ and $\displaystyle\lim_{k\to\infty}u_k=\infty$.
 If $u$ is a positive solution of $(\ref{A1})$ in $Q_{\infty}$, it possesses an initial trace which is  a positive Radon measure $\gm$ on $\BBR^N$. This result holds in particular if $(C1)$ and $(C3)$ hold.
\es

The proofs are combination of methods developed in \cite{MV4} for elliptic equations, stability results and \rth{ALA} and \rth{TTA}. 
\mysection{Isolated singularities}

In order to study $(\ref{A1})$, we start proving \rprop{exist-fundamental}. \medskip

\noindent{\bf Proof of \rprop{exist-fundamental}} \smallskip

We denote by $E(x,t)=(4\pi t)^{-N/2}e^{-\abs{x}^2/4t}$ the fundamental solution of the heat equation in $Q_\ity$. Since $kE$  ($k>0$) is a supersolution for $(\ref{A1})$, it is classical to prove that if 
\bel{I1} I:=\int_0^1\int_{B_R}f(kE(x,t))dx\,dt < \ity \ee
for any $R>0$, then there exists a unique solution $u=u_k$ to $(\ref{A1})$ satisfying initial condition 
$u_k(.,0)=k\gd_0 $
in $\CD'(\BBR^N)$. 
Furthermore the mapping $k \mapsto u_k$ is increasing. Actually, it is proved in \cite[Th 1.1]{MV2} that if $f$ satisfies the weak singularity assumption $(\ref{A8})$, then for any positive bounded Borel measure there exists a unique solution $u:=u_\gm$ to $\ref{A1}$ satisfying $u_\gm(.,0)=\gm$. Furthermore if $\{\gm_n\}$ is a sequence of positive bounded measures which converge to a measure $\gm$ in the weak-star topology of measures, then the sequence of corresponding solutions $\{u_{\gm_n}\}$ converges locally uniformly to $u_\gm$, and $\{f(u_{\gm_n})\}$ converges to $f(u_\gm)$ in $L^1_{loc}(\BBR^N\ti [0,\infty))$.

This existence result and the next proposition lead to the conclusion of \rprop{exist-fundamental}. \qeda
\bprop{E-I} If $f$ satisfies $(\ref{A8})$ and $(C1)$ then $(\ref{I1})$ is fulfilled.
\es
\Proof We set 
\bel{h} h(r)=\myfrac{f(r)}{r} \qq r \in (0,\ity).\ee
$I$ is rewritten as
\[ I=kC^*\int_0^1\int_{B_R}t^{-N/2}e^{-{\abs x}^2/4t}h(kC^*t^{-N/2}e^{-{\abs x}^2/4t})dx\,dt \]
where $C^*=(4\pi)^{-N/2}$. Put $r=\abs x$ then $dx=r^{N-1}dr$, and 
\[ I=kC^*\int_0^1t^{-N/2}\int_0^R e^{-r^2/4t}h(kC^*t^{-N/2}e^{-r^2/4t})r^{N-1}dr\,dt. \]
Wet put $\gr=\frac{r}{\sqrt{t}}$, then $r^{N-1}dr=\gr^{N-1}t^{N/2}d\gr$, and
\[ I=kC^*\int_0^1\int_0^{R/\sqrt{t}}e^{-\gr^2/4}h(kC^*t^{-N/2}e^{-\gr^2/4})\gr^{N-1}d\gr \, dt. \]
We set
\[ I_1:=kC^*\int_0^1\int_0^1e^{-\gr^2/4}h(kC^*t^{-N/2}e^{-\gr^2/4})\gr^{N-1}d\gr \,dt, \]	 
\[ I_2:=kC^*\int_0^1\int_1^{R/\sqrt{t}} e^{-\gr^2/4}h(kC^*t^{-N/2}e^{-\gr^2/4})\gr^{N-1}d\gr \,dt. \]
Since $e^{-\gr^2/4}\gr^{N-1}$ is bounded in $[0,\ity)$, then there exists a constant $c_1$ depending only on $k$ such that
\[ I_1 < c_1 \int_0^1\int_0^1 h(kC^*t^{-N/2})d\gr \,dt = c_1\int_0^1h(kC^*t^{-N/2})dt < \ity. \]
Next we show that under the condition $(\ref{A8})$, $I_2 < \ity$. In order to do that we introduce the variable $\gt$ such that $t^{-N/2}e^{-\gr^2/4}=\gt^{-N/2}$. Then $t=\gt e^{-\frac{\gr^2}{2N}}$ and $dt=e^{-\frac{\gr^2}{2N}}d\gt$. Therefore 
\bel{E-I1} I_2 \leq kC^*\int_1^\ity e^{-\frac{(N+2)\gr^2}{4N}}\gr^{N-1}\bigg(\int_0^{e^{\gr^2/2N}}h(kC^*\gt^{-N/2})d\gt\bigg) d\gr. 
\ee
Since $h$ satisfies $(\ref{A8})$, there exists $\ge >0$ (depending only on $k$) such that  
\[ \int_0^\ge h(kC^*\gt^{-N/2})d\gt \]
take a finite value, denoted by $c_2$. Hence
\bel{E-I2} \int_0^{e^{\gr^2/2N}}h(kC^*\gt^{-N/2})d\gt \leq c_2 + h(kC^*\ge^{-N/2})(e^\frac{\gr^2}{2N}-\ge). 
\ee
Inserting $(\ref{E-I2})$ into the right-hand side of $(\ref{E-I1})$, we obtain
\[ I_2 \leq c_3\int_1^\ity e^{-\frac{(N+2)\gr^2}{4N}}\gr^{N-1} d\gr + c_4 \int_1^\ity e^{-\frac{\gr^2}{4}}\gr^{N-1} d\gr < \ity \]
where $c_3=kC^*c_2$ and $c_4=kC^*h(kC^*\ge^{-N/2})$. Thus $I=I_1+I_2<\ity$. \qeda \medskip

The functions which satisfy the following ODE are particular solutions of $(\ref{A1})$
	\bel{ODE}
		 \prt_t\gf+f(\gf) = 0 \qq \text{in } (0,\ity).
	\ee
For $a>0$, we denote by $\gf_a$ the solution of $(\ref{ODE})$ with initial data $\gf(0)=a$. If $(\ref{A15})$ holds then ${\displaystyle \lim_{a \to \ity}}\gf_a(t)=\ity$ for any $t \in (0,\ity)$. While, if $(\ref{A10})$ holds
there exists a maximal solution $\gf_\ity$ given explicitely by
\[ t=\int_{\gf_\ity(t)}^\ity\frac{ds}{f(s)}<\ity. \]	
\blemma{J} If $(\ref{A15})$ holds then 
$${\displaystyle \liminf_{r \to \infty}} \frac{f(r)}{r\ln^\ga r} =0,  \forevery \ga>1. $$ 
If $(\ref{A10})$ holds then 
$${\displaystyle \limsup_{r \to \infty}}\,\frac{f(r)}{r\ln^\ga r} =\ity, \forevery 0<\ga \leq 1.$$
\es
\Proof \noindent{\it Case 1.} Assume $(\ref{A15})$ holds then
\bel{J} J:=\myint{e}{\ity}\frac{ds}{f(s)} < \ity. \ee 
We put $s=e^{r^{-1}}$ and derive 
\[ J=\int_0^1\frac{dr}{r^2h(e^{r^{-1}})} \]
where $h$ is defined in $(\ref{h})$. Suppose that there exists $\ga >1$ such that 
\[ \liminf_{s \to \infty}\myfrac{f(s)}{s\ln^\ga s}> 0, \]
equivalently,
\[ \liminf_{r \to 0}r^\ga h(e^{r^{-1}})>0, \]
then there exists $l>0$ and $r_0 \in (0,1)$ such that
\[ h(e^{r^{-1}}) > lr^{-\ga} \qq \forall r \in (0,r_0). \] 
Hence we derive the following contradiction
\[ J<\frac{1}{l}\int_0^{r_0}r^{\ga-2}dr + \int_{r_0}^1\frac{dr}{r^2h(e^{r^{-1}})} < \ity. \]

\noindent{\it Case 2.} Assume $(\ref{A10})$ holds then $J=\ity$. Suppose that there exists $\ga \in (0,1]$ such that 
\[ \limsup_{s \to \infty}\myfrac{f(s)}{s\ln^\ga s}< \ity, \]
equivalently,
\[ \limsup_{r \to 0}r^\ga h(e^{r^{-1}})<\ity, \]
then there exists $l>0$ and $r_0 \in (0,1)$ such that
\[ h(e^{r^{-1}}) < lr^{-\ga} \qq \forall r \in (0,r_0). \] 
Hence
\[ J>\frac{1}{l}\int_0^{r_0}r^{\ga-2}dr + \int_{r_0}^1\frac{dr}{r^2 h(e^{r^{-1}})} = \ity, \]
which is a contradiction.\qeda \medskip

\noindent{\bf Proof of \rth{ALA}.} \medskip

Since $(\ref{A15})$ holds, by \rlemma{J} and the definition $(\ref{h})$ of $h$, 
$${\displaystyle \liminf_{r \to \ity}}\myfrac{h(r)}{\ln^\ga r}=0 \forevery \ga>1.$$ 
Thus
\[ \liminf_{r \to \ity}\myfrac{h(r)}{\ln^\ga r}=0 \forevery \ga>2. \]
By (C3), there exists $\gb \in (1,2]$ such that ${\displaystyle \limsup_{r \to \ity}}\,h(r)/\ln^\gb r < \ity$. Hence there exist $M>0$ and $r_0>0$ such that
\bel{ALA-1} h(r) < M\ln^\gb r \qq \forall r \in (r_0,\ity). \ee 

\noindent{\it Step 1.} Let $k >0$, we claim that
\bel{ALA-1'}  \gth_k(t)<2^{\gb-1}Mt(\ln  k)^\gb+\frac{MN^\gb}{2}\int_0^1(\ln  (\gt^{-1}))^\gb d\gt \qq \forall t \in (0,1)  \ee
where $\gth_k(t)=\myint{0}{t} h(kC^*\gt^{-N/2})d\gt$ with $C^*=(4\gp)^{-N/2}$.
Set $r=kC^*\gt^{-N/2}$ then $(\ref{ALA-1})$ becomes 
\[ h(kC^*\gt^{-N/2}) < M[\ln (kC^*)+\frac{N}{2}\ln (\gt^{-1})]^\gb \qq \forall \gt \in (0,\gt_0) \]
where $\gt_0=(kC^*)^{2/N}r_0^{-2/N}$. 
We put $a_1=\ln k$, $a_2=\frac{N}{2}\ln(\gt^{-1})$, and apply the following inequality
\[ (a_1+a_2)^\gb \leq 2^{\gb-1}(a_1^\gb+a_2^\gb) \]
in order to obtain
\bel{ALA-2} \BA{l} h(kC^*\gt^{-N/2}) < M[\ln (k)+\frac{N}{2}\ln (\gt^{-1})]^\gb \\[2mm]
\phantom{ h(kC^*\gt^{-N/2})}
\leq 2^{\gb-1}M[(\ln k)^\gb+(\frac{N}{2})^\gb \ln ^\gb(\gt^{-1})] \qq \forall \gt \in (0,\gt_0),
\EA \ee
 (notice that $C*=(4\pi)^{-N/2}<1$).  Integrating over $[0,t]$ yields to $(\ref{ALA-1'})$. \smallskip

\noindent{\it Step 2.} It follows from $(\ref{ALA-2})$ that $(\ref{A8})$ is fulfilled; hence by \rprop{exist-fundamental} there exists a unique solution $u_k$ of $(\ref{A1})$ in $Q_\ity$ with initial data $k\gd_0$. By the maximum principle, $u_k(x,t) \leq kE(x,t)$ for every $(x,t) \in Q_\ity$, which implies $u_k(x,t) \leq kC^*t^{-N/2}$ for every $(x,t) \in Q_\ity$. Therefore, since $h$ is increasing, 
\[ \prt_t u_k - \Gd u_k + u_k h(kC^*t^{-N/2}) \geq 0. \]
If we set $v_k(x,t)=e^{\gth_k(t)}u_k(x,t)$, we obtain
\[ \prt_t v_k -\Gd v_k = e^{\gth_k(t)}[\prt_t u_k - \Gd u_k + u_k h(kC^*t^{-N/2})] \geq 0 \]
and $v_k(.,0) = u_k(.,0)=k\gd_0$. By the maximum principle, there holds
\bel{ALA-3} v_k(x,t) \geq kC^*t^{-N/2} e^{-{\abs x}^2/4t} \Llra u_k(x,t) \geq kC^*t^{-N/2}e^{-\gth_k(t)-{\abs x}^2/4t}. \ee
By step 1,
\bel{ALA-4} e^{-\gth_k(t)} \geq c_1e^{-M_\gb t(\ln k)^\gb} \qq \forall t \in (0,1) \ee
where 
$$ c_1=\exp \Big(-\frac{M(N)^\gb}{2}\int_0^1(\ln (\gt^{-1}))^\gb d\gt \Big)$$ 
and $M_\gb=M2^{\gb-1}$. Inserting $(\ref{ALA-4})$ into the right-hand side of $(\ref{ALA-3})$, we get
\[ u_k(x,t) \geq c_1C^*t^{-N/2}e^{\ln k-M_\gb t(\ln k)^\gb-{\abs x}^2/4t} \qq \forall (x,t) \in Q_1:=\BBR^N \ti (0,1).\]
If ${\dsps \lim_{k \to \ity}}u_k(x,t)<\ity$ for all $(x,t) \in Q_\ity$, we put $\unl U:={\dsps \lim_{k \to \ity}u_k}$, then
\[ \unl U(x,t) \geq c_1C^*t^{-N/2}e^{\ln k-M_\gb t(\ln k)^\gb-{\abs x}^2/4t}\q \forall (x,t) \in Q_1,\q \forall k>0.\]
Let $\{t_n\} \sbs (0,1]$ be a sequence converging to $0$. We choose $k_n=\exp\left((2M_\gb t_n)^\frac{1}{1-\gb}\right)$
then $\ln k_n-M_\gb t_n(\ln k_n)^\gb = \frac{1}{2}\ln k_n$. Next we restrict $x$ in order to have
	\[ \ln k_n-M_\gb t_n(\ln k_n)^\gb-\frac{{\abs x}^2}{4t_n}=\frac{1}{2}\ln k_n-\frac{{\abs x}^2}{4t_n} \\ \geq 0 \Llra \abs x \leq 2^\frac{\gb-2}{2(\gb-1)}M_\gb^\frac{1}{2(1-\gb)}t_n^\frac{\gb-2}{2(\gb-1)}. \]
Therefore, since $1< \gb \leq 2$,
	\[ \lim_{n \to \ity}\unl U (x,t_n) =\ity \]
uniformly on ${\BBR}^N$ if $1\leq \gb <2$, or uniformly on the ball $B_{r_2}$ where $r_2=(2M)^{-1/2}$ if $\gb=2$. Since the sequence $\{t_n\}$ is arbitrary, 
	\[ \lim_{t \to 0}\unl U (x,t) =\ity \] 
uniformly on ${\BBR}^N$ if $1\leq \gb <2$, or uniformly on the ball $B_{r_2}$ if $\gb=2$. \smallskip

We pick some point $x_0$ in ${\BBR}^N$ (resp. $B_{r_2}$) if $1<\gb<2$ (resp. $\gb=2$). Since for any $k>0$, the solution $u_{k\gd_{x_0}}$ of $(\ref{A1})$ with initial data $k\gd_{x_0}$ can be approximated by solutions with bounded initial data and support in $B_\gs (x_0)$ where $0<\gs<r_2 - \abs {x_0}$, it follows
	\[ \unl U(x,t) \geq u_{k\gd_{x_0}}(x,t) = u_k(x-x_0,t),  \]
by comparison principle. Letting $k \to \ity$ yields to $\unl U(x,t) \geq \unl U(x-x_0,t)$. Interverting the role of $0$ and $x_0$ yields to $\unl U(x,t) = \unl U(x-x_0,t)$. If we iterate this process we derive
	\[ \unl U(x,t) = \unl U(x-y,t) \forevery y \in {\BBR}^N. \]
This implies that $\unl U(x,t)$ is independent of $x$ and therefore it is a solution of $(\ref{A11})$. By $(\ref{A15})$, $\unl U(x,t) = \ity$ for any $(x,t) \in Q_\ity$, which is a contradiction and the conclusion follows. \qeda 
\bprop{Up1} Assume $(\ref{A10})$ is satisfied. For any $k>0$, there holds
\[ u_k(x,t) \leq \gf_\ity(t) \qq \forall (x,t) \in Q_\ity. \]
\es
\Proof For any small $\ge >0$, we set $\gf_{\ity\ge}(t)=\gf_\ity(t-\ge), t \in \left[\ge,\ity\right)$ then $\gf_{\ity\ge}$ is a solution of $(\ref{A1})$ in $(\ge,\ity)$, which dominates $u_k$ on ${\BBR}^N \ti \{\ge\}$ for any $k>0$. By comparison principle, $u_k(x,t) \leq \gf_{\ity\ge}(t)$ for every $(x,t) \in \BBR^N \ti \left[\ge,\ity\right)$. Letting $\ge \to 0$ yields the claim. \qeda 

A necessary and sufficient condition for the existence of a maximal solution to the stationary equation
\[ -\Gd w +f(w) = 0 \]
in a bounded domain $\Gw$ is the Keller-Osserman condition $(\ref{A12})$ (\cite{Ke}, \cite{Os}). If $f$ is convex and $(\ref{A12})$ holds, then $(\ref{A10})$ is fulfilled. The Keller-Osserman condition can be replaced by another condition, which owes to the following result.
\blemma{KON} Assume $f$ is convex on $(0,\ity)$. Set
\[ L:=\int_1^\ity\frac{ds}{\sqrt{sf(s)}}. \]
Then $(\ref{A12})$ holds if and only if $L < \ity$.
\es
\Proof In order to obtain the assertion, it is sufficient to show that 
\bel{KON1} s\,f(\frac{s}{2}) \leq F(s) \leq s\,f(s) \qq \forall s\geq 1. \ee
The right-hand side estimate in $(\ref{KON1})$ follows from the monotone property of $f$. The assumption of convexity of $f$ in $(0,\ity)$ implies
\[ f(s) \geq f(\frac{s}{2})+\frac{s}{2}f\pr(\frac{s}{2}) \qq \forall s>0. \]
Define $\vgf(s)=\myint{0}{s}f(\gs)d\gs-sf(\frac{s}{2})$, then $\vgf\pr(s)=f(s)-f(\frac{s}{2})-\frac{s}{2}f\pr(\frac{s}{2}) \geq 0$. Hence $\vgf(s) > \vgf(0)=0$, which leads to the left-hand side estimate in $(\ref{KON1})$. \qeda \medskip

By using the same argument as in the proof of the \rlemma{J} and thank to the \rlemma{KON}, we obtain the following lemma.
\blemma{K} If $(\ref{A16})$ holds then 
$${\displaystyle \liminf_{r \to \infty}} \myfrac{f(r)}{r\ln^\ga(r)}=0 \forevery \ga>2.$$
If $(\ref{A12})$ holds then 
$${\displaystyle \limsup_{r \to 0}}\myfrac{f(r)}{r\ln^\ga(r)}=\ity \forevery 0<\ga \leq 2.$$
\es
\noindent {\bf Proof of \rth{TTA}.} \medskip

Since $(\ref{A16})$ holds, by \rlemma{K} and the definition $(\ref{h})$ of $h$, 
$${\displaystyle \liminf_{r \to \ity}}{h(r)}{\ln^\ga r}=0 \forevery \ga>2. 
$$
By (C3), there exists $\gb \in (1,2]$ such that ${\displaystyle \limsup_{r \to \ity}}\,h(r)/\ln^\gb r < \ity$. Hence there exists $M>0$ and $r_0>0$ such that
\bel{TTA-1} h(r) < M\ln^\gb r \forevery r \in (r_0,\ity). \ee 
\noindent{\it Step 1.} For any $k >0$ we set
$$\gth_k(t)=\int_0^t h(kC^*\gt^{-N/2})d\gt
$$
where $C^*=(4\gp)^{-N/2}$. We claim that
\bel{TTA-1'}  
\gth_k(t)<2^{\gb-1}Mt(\ln k)^\gb+\frac{MN^\gb}{2}\int_0^1(\ln (\gt^{-1}))^\gb d\gt \q \forall t \in (0,1). 
\ee 
If we define $\gt$ by $r=kC^*\gt^{-N/2}$, $(\ref{TTA-1})$ becomes 
\[ h(kC^*\gt^{-\frac{N}{2}}) < M[\ln (kC^*)+\frac{N}{2}\ln (\gt^{-1})]^\gb \q \forall \gt \in (0,\gt_0) \]
where $\gt_0=(kC^*)^{2/N}r_0^{-2/N}$. 
We set $a_1=\ln k$, $a_2=\frac{N}{2}\ln (\gt^{-1})$, and apply the following inequality
\[ (a_1+a_2)^\gb \leq 2^{\gb-1}(a_1^\gb+a_2^\gb) \]
in order to obtain (notice that $C^*<1$)
\bel{TTA-2} \BA{l} h(kC^*\gt^{-N/2}) < M[\ln (k)+\frac{N}{2}\ln (\gt^{-1})]^\gb \\
\phantom{ h(kC^*\gt^{-N/2})}
\leq 2^{\gb-1}M[(\ln k)^\gb+(\frac{N}{2})^\gb \ln ^\gb(\gt^{-1})].
\EA \ee
Integrating over $[0,t]$, we obtain $(\ref{TTA-1'})$. \smallskip

\noindent{\it Step 2.} It follows from $(\ref{TTA-2})$ that $(\ref{A8})$ is fulfilled; hence by \rprop{exist-fundamental} there exists a unique solution of $(\ref{A1})$ in $Q_\ity$ with initial trace $k\gd_0$. By maximum principle, $u_k(x,t) \leq kE(x,t)$ for every $(x,t) \in Q_\ity$, which implies that $u_k(x,t) \leq kC^*t^{-N/2}$ for every $(x,t) \in Q_\ity$. Therefore, since $h$ is increasing, 
\[ \prt_t u_k - \Gd u_k + u_k h(kC^*t^{-N/2}) \geq 0. \]
We set $v_k(x,t)=e^{\gth_k(t)}u_k(x,t)$ and obtain
\[ \prt_t v_k -\Gd v_k = e^{\gth_k(t)}[\prt_t u_k - \Gd u_k + u_k h(kC^*t^{-N/2})] \geq 0, \]
with $v_k(.,0) = u_k(.,0)=k\gd_0$. By maximum principle, it follows
\bel{TTA-3} v_k(x,t) \geq kC^*t^{-N/2} e^{-{\abs x}^2/4t} \Llra u_k(x,t) \geq kC^*t^{-N/2}e^{-\gth_k(t)-{\abs x}^2/4t}. \ee
By step 1,
\bel{TTA-4} e^{-\gth_k(t)} \geq c_1e^{-M_\gb t(\ln k)^\gb} \q \forall t \in (0,1) \ee
where 
$ c_1=\exp \Big(-\frac{M(N)^\gb}{2}\int_0^1(\ln (\gt^{-1}))^\gb d\gt \Big) $ and $M_\gb=M2^{\gb-1}$. 
Inserting $(\ref{TTA-4})$ into the right-hand side of $(\ref{TTA-3})$, we get
\[ u_k(x,t) \geq c_1C^*t^{-N/2}e^{\ln k-M_\gb t(\ln k)^\gb-{\abs x}^2/4t} \forevery (x,t) \in Q_1=\BBR^N \ti (0,1) .\]
Since $k \mapsto u_k$ is increasing, by \rprop{Up1} there exists $\unl U:={\displaystyle \lim_{k \to \ity}u_k}$ and  $\unl U \geq u_k$. Hence
\[ \unl U(x,t) \geq c_1C^*t^{-N/2}e^{\ln k-M_\gb t(\ln k)^\gb-{\abs x}^2/4t} \forevery (x,t) \in Q_1, \forall k>0. \]
Let $\{t_n\} \sbs (0,1]$ be a sequence converging to $0$. We choose $k_n=\exp ((2M_\gb t_n)^\frac{1}{1-\gb})$, equivalently $\ln k_n-M_\gb t_n(\ln k_n)^\gb = \frac{1}{2}\ln k_n$. Next we restrict $\abs x$ in order 
	\[ \ln k_n-M_\gb t_n(\ln k_n)^\gb-\frac{{\abs x}^2}{4t_n}=\frac{1}{2}\ln k_n-\frac{{\abs x}^2}{4t_n} \geq 0 \Llra \abs x \leq r_\gb\, t_n^\frac{\gb-2}{2(\gb-1)}, \]
where $r_\gb=2^\frac{\gb-2}{2(\gb-1)}M_\gb^\frac{1}{2(1-\gb)}$. 
Because $1< \gb \leq 2$, it follows
	\[ \lim_{n \to \ity}\unl U (x,t_n) =\ity, \]
uniformly on ${\BBR}^N$ if $1\leq \gb <2$, or uniformly on the ball $B_{r_2}$ where $r_2=(2M)^{-\frac{1}{2}}$ if $\gb=2$. Since the sequence $\{t_n\}$ is arbitrary, 
	\[ \lim_{t \to 0}\unl U (x,t) =\ity \] 
uniformly on ${\BBR}^N$ if $1\leq \gb <2$, or uniformly on the ball $B_{r_2}$ if $\gb=2$. \smallskip

We pick some point $x_0$ in ${\BBR}^N$ (resp. $B_{r_2}$) if $1<\gb<2$ (resp. $\gb=2$). Since for any $k>0$, the solution $u_{k\gd_{x_0}}$ of $(\ref{A1})$ with initial data $k\gd_{x_0}$ can be approximated by solutions with bounded initial data and support in $B_\gs (x_0)$ where $0<\gs<r_2 - \abs {x_0}$,  it follows
	\[ \unl U(x,t) \geq u_{k\gd_{x_0}}(x,t) = u_k(x-x_0,t), \]
by comparison principle. Letting $k \to \ity$ yields to $\unl U(x,t) \geq \unl U(x-x_0,t)$. Reversing the role of $0$ and $x_0$ yields to $\unl U(x,t) = \unl U(x-x_0,t)$. If we iterate this process we derive
	\[ \unl U(x,t) = \unl U(x-y,t) \forevery y \in {\BBR}^N. \]
This implies that $\unl U(x,t)$ is independent of $x$ and therefore it is a solution of $(\ref{A11})$
Since $(\ref{A10})$ holds, $\unl U(x,t) = \gf_\ity(t)$ for every $(x,t) \in Q_\ity$. \qeda 

\bprop{Up2} Assume $(\ref{A12})$ and $(\ref{A8})$ are satisfied. Then for any $k>0$ there holds
\[  u_k(x,t) \leq \Gf(\abs x) \forevery (x,t) \in Q_\ity \]
where $\Gf$ is a solution to the problem
\[ \left\{ \BA{ll}
-\Gf \ppr + f(\Gf) &= 0 \qq \text{ in } (0,\ity) \\\phantom{--,}
{\displaystyle \lim_{s \to 0}}\Gf (s)&=\ity.
\EA \right. \]
\es 
\Proof \noindent {\it Step 1: Upper estimate.} Since $f$ satisfies $(\ref{A12})$, by \cite{Ke} for any $R>0$, there exists a solution $w_R$ to the problem
\bel{w_R} \left\{ \BA{ll}
- \Gd w_R + f(w_R) &= 0 \qq \text{ in } B_R, \\\phantom{-;;;}
{\displaystyle \lim_{\abs x \to R}}w_R(x)&=\ity,
\EA \right. \ee
and $w_R$ is nonnegative since $f(0)=0$. Notice also that $R\mapsto w_R$ is decreasing, since $f$ is nondecreasing; moreover $\lim_{R\to\infty}w_R=0$, since $f(0)=0$ and $f$ is positive on $(0,\infty)$. 
Let $x_0 \ne 0$ arbitrary in ${\BBR}^N$. Set $\BBE=\{\vec{e}: \abs{\vec{e}}=1\}$ and take $\vec{e} \in \BBE $. Put $x_{\vec{e}}=\abs {x_0}\vec{e}$ and for $n > \abs {x_0}$ put $a_n=n\vec{e}$. Denote by ${\BBH}_{\vec{e}}$ the open half-space generated by $\vec{e}$ and its orthogonal hyperplane at the origin, then $x_{\vec{e}},a_n \in {\BBH}_{\vec{e}}$. Take $R$ such that $n-\abs {x_0} < R < n$. We set $W_{\vec{e},n,R}(x)=w_R(x-a_n)$, then $W_{\vec{e},n,R}$ is a solution of $(\ref{A1})$ in $B_R(a_n)$ and blows-up on the boundary ${\displaystyle \lim_{\abs {x-a_n} \to R}}W_{\vec{e},n,R}(x)=\ity$. By the maximum principle,  
\bel{UB1} u_k(x,t) \leq W_{\vec{e},n,R}(x) \forevery (x,t) \in B_R(a_n) \ti (0,\ity). 
\ee 
The sequence $\{W_{\vec{e},n,R}\}$ is decreasing with respect to $R$ and is bounded from below by $u_k$, then there exists $W_{\vec{e},n}:={\displaystyle \lim_{R \to n}}W_{\vec{e},n,R}$ satisfying
\bel{UB2} u_k(x,t) \leq W_{\vec{e},n}(x) \forevery (x,t) \in B_n(a_n)\ti (0,\ity). 
\ee 
The sequence $\{W_{\vec{e},n}\}$ is also decreasing with respect to $n$ and is bounded from below by $u_k$, then there exists $W_{\vec{e},\ity}:={\displaystyle \lim_{n \to \ity}}W_{\vec{e},n}$. Letting $n \to \ity$ in $(\ref{UB2})$ yields to
\bel{UB3} u_k(x,t) \leq W_{\vec{e},\ity}(x) \q \forall (x,t) \in {\BBH}_{\vec{e}}\ti (0,\ity). \ee  
In particular, 
\bel{UB4} u_k(x_{\vec{e}},t) \leq W_{\vec{e},\ity}(x_{\vec{e}}). \ee
Since $u_k$ is radial, it follows that 
\[ u_k(x_0,t)=u_k(x_{\vec{e}},t) \leq W_{\vec{e},\ity}(x_{\vec{e}}). \]
For any $r>0$, $n>r$, $n-r<R<n$ and $\vec e$, $\vec{e'} \in \BBE$, since $w_R$ is radial,
\[ w_R(r\vec{e}-n\vec{e})=w_R(r\vec{e'}-n\vec{e'}). \]
Letting successively $R \to n$, $n \to \ity$ yields to
\[  W_{\vec{e},\ity}(r\vec{e})=W_{\vec{e'},\ity}(r\vec{e'}). \]
Define $\tl \Gf(r):=W_{\vec{e},\ity}(r\vec{e})$, $\forall r \in (0,\ity)$ then it satisfies
\bel{Phi-til} \left\{ \BA{ll}
-\tl \Gf\ppr -\myfrac{N-1}{r} \tl \Gf \pr + f(\tl \Gf) &= 0 \qq \text{in } (0,\ity) \\ \phantom{-------}
{\displaystyle \lim_{r \to 0}}\tl \Gf(r)&=\ity ,
\EA \right.\ee  
and 
\bel{UB5} u_k(x,t) \leq \tl \Gf(\abs x) \qq \forall (x,t)\in Q_\ity.\ee
\noindent {\it Step 2: End of the proof.}  We claim that 
\bel{UB6} \tl \Gf(r)\leq \Gf(r) \qq \forall r \in (0,\ity). \ee
For any $\ge>0$, we set $\Gf_\ge(r)=\Gf(r-\ge)$, $r>\ge$ then $\Gf_\ge$ is a solution of 
\bel{UB7} -\Gf\ppr_\ge + f(\Gf_\ge) = 0 \qq \text{in } (\ge,\ity) \ee
verifying ${\displaystyle \lim_{r \to \ge}}\Gf_\ge(r)=\ity$. Since $\Gf\pr_\ge \leq 0$, $\Gf_\ge$ is a supersolution of the equation in $(\ref{Phi-til})$ in $(\ge,\ity)$, which dominates $\tl \Gf$ at $r=\ge$. By the maximum principle, $\tl \Gf \leq \Gf_\ge$ in $(\ge,\ity)$. Letting $\ge \to 0$ yields $(\ref{UB6})$. Combining $(\ref{UB5})$ and $(\ref{UB6})$ leads to the conclusion. \qeda \medskip

\noindent {\it Remark.} Combining \rprop{Up1} and \rprop{Up2} yields to
\bel{UNIV-k} u_k(x,t) \leq min\{\gf_\ity(t),\Gf(\abs x)\} \qq \forall (x,t) \in Q_\ity, \forall k>0. 
\ee  \smallskip

\noindent{\bf Proof of \rth{min-solution}.} \smallskip

Since $f$ is convex, (\ref{A12}) implies (\ref{A10}). Actually, only $\liminf_{s\to\infty}\frac{f(s)}{s}>0$ is needed for this implication. The sequence $\{u_k\}$ is increasing with respect to $k$ and is bounded from above by $(\ref{UNIV-k})$ then there exists $\unl U:={\dsps \lim_{k \to \ity}u_k}$ satisfying
\bel{UNIV} \unl U(x,t) \leq min\{\gf_\ity(t),\Gf(\abs x)\} \qq \forall (x,t) \in Q_\ity, \forall k>0. 
\ee
Moreover, $\unl U \in \CU_0$ because $\unl U$ has the following properties:\\
\noindent (i) It is positive in $Q_\ity$, belongs to $C(\ovl Q \sms\{(0,0)\})$ and vanishes on ${\BBR}^N \ti \{(0,0)\} \sms \{0\}$.\\
\noindent (ii) It satisfies $(\ref{A1})$ and 
\bel{UNIV1} \lim_{t \to 0}\int_{B_\gs}\unl U(x,t)dx=\ity, \q \forall \gs > 0. \ee 
In the sense of initial trace in \rdef{Trace}, $\unl U$ has initial trace $tr_{\BBR^N}(\unl U)=(\{0\},0)$ (here $\{0\}$ is the singular part and the Radon measure on ${\BBR}^N \sms \{0\}$ is the zero measure) and the conclusion follows from \rprop{singular-point}. \qeda \medskip

By a simple adaptation of the proof of \rprop{Up1} and \rprop{Up2} it is possible to extend $(\ref{UNIV})$ to any positive solution vanishing on ${\BBR}^N \ti \{0\} \sms \{(0,0)\}$.

\bprop{Up2*} Assume $(\ref{A12})$ and $(C2)$ are satisfied. Then any positive solution $u\in C^{2,1}(Q_\infty)$ of $(\ref{A1})$ satisfies
\bel{UNIVE} u(x,t) \leq \gf_\ity(t) \qq \forall (x,t) \in Q_\ity. \ee
If we assume moreover that  $u\in C(\ovl Q \sms\{(0,0)\})$ vanishes on ${\BBR}^N \ti \{0\} \sms \{0\}$, there holds
\bel{UNIV'} u(x,t) \leq min\{\gf_\ity(t),\Gf(\abs x)\} \qq \forall (x,t) \in Q_\ity. \ee
\es 
\Proof Since $f(0)=0$ and due to the convexity of $f$, the following inequality holds
\bel{convex-ineq}
f(a+b) \geq f(a)+f(b) \forevery a,b >0,
\ee
which implies that for any $R,\gt>0$, $(x,t)\mapsto \phi_\infty(t-\gt)+w_R(x)$ is a supersolution of $(\ref{A1})$ in $B_R\ti (\gt,\infty)$. This function dominates $u$ on the parabolic boundary, thus in the domain itself by the comparison principle. Since $f(r)>0$ if $r>0$, ${\displaystyle \lim_{R\to\infty}} w_R=0$ in $\BBR^N$. Therefore
$$u(x,t)\leq \phi_\infty(t)=\lim_{\gt\to 0}\lim_{R\to\infty}(\phi_\infty(t-\gt)+w_R(x))\qq\forall (x,t)\in Q_\infty.
$$
For the second estimate we notice that $(\ref{UB1})$ is valid with $u_k$ replaced by $u$ (and without assumption $(\ref{A8})$ since existence is assumed). The remaining of the proof of \rprop{Up2} is similar and yields to 
$$u(x,t)\leq \Gf(|x|)\qq\forall (x,t)\in Q_\infty.
$$
This implies $(\ref{UNIV'})$.\qeda

\medskip

It is also possible to construct a maximal element of $\CU_0$ ($\CU_0$ is defined in \rth{min-solution}). For $\ell >0$ and $\ge>0$, let $u:=U_{\ge,\ell}$ be the solution of 
\[ \left\{ \BA{lll} 
\prt_t u - \Gd u + f(u) &= 0 \qq &\text{in } Q_\ity\\
\phantom{\prt_t u - \Gd u +}
u(x,0)&=\ell\gc_{B_\ge} &\text{in } {\BBR}^N.
\EA \right. \]
\blemma{Comparison1} For any $\gt>0$ and $\ge>0$, there exist $\ell >0$ and $m(\gt,\ge)>0$ such that any positive solution $u$ of $(\ref{A1})$ which verifies (i) in the proof of \rth{min-solution} satisfies
	\bel{Comparison1.1} u(x,t)\leq U_{\ge,\ell}(x,t-\gt)+m(\gt,\ge) \qq \forall (x,t)\in Q_\ity,\, t\geq\gt. \ee
Furthermore 
	\be \lim_{\gt \to 0}m(\gt,\ge)=0 \qq \forall\ge>0.	\ee
Finally 
\bel{maxs}
\overline U(x,t)=\lim_{\gt\to 0}\lim_{\ge\to 0}\lim_{\ell\to \ity}\left(U_{\ge,\ell}(x,t-\gt)+m(\gt,\ge)\right)\ee
is the maximal element of $\CU_0$.
\es
\Proof We set $\ell=\gf_\ity(\gt)$, then $u(x,\gt)\leq\ell$ for any $x\in{\BBR}^N$. Let $W:=W_{\ge/2}$ be the solution of the following Cauchy-Dirichlet problem
\bel{eW} \left\{ \BA{lll}
\prt_tW-\Delta W+f(W) &= 0 \qq &\text{in } B_{\ge/2}^c\ti(0,\ity)\\
\phantom{\prt_tW-\Delta -,,}
W(x,0)&=0  &\text{in } B_{\ge/2}^c\\
\phantom{\prt_tW-\Delta -,,}
W(x,t)&=\gf_\ity(t) \q &\text{in } \prt B_{\ge/2}^c\ti(0,\ity)
\EA \right. \ee
and put $ m(\gt,\ge):=max\{W_{\ge/2}(x,\gd):\abs x >\ge, 0<\gd\leq\gt\}$. It is clear to see that 
	\be \lim_{\gt \to 0}m(\gt,\ge)=W_{\ge/2}(x,0)=0. \ee
From the fact that $u(x,0)=0 \text{ in } B_{\ge/2}^c$, $u(x,t)\leq\gf_\ity(t) \text{ in } \prt B_{\ge/2}^c\ti(0,\ity)$ and the maximum principle, it follows that $u(x,t)\leq W_{\ge/2}(x,t) \text{ in } B_{\ge/2}^c \ti (0,\ity)$. 

Next, we compare $U_{\ge,\ell}(.,.-\gt)+m(\gt,\ge)$ with $u$ in ${\BBR}^N \ti (\gt,\ity)$. The function $U_{\ge,\ell}(.,.-\gt)+m(\gt,\ge)$ is a supersolution of $(\ref{A1})$ in ${\BBR}^N \ti (\gt,\ity)$. If $x \in B_\ge$,  $U_{\ge,\ell}(x,0)=\ell\geq u(x,\gt)$, which implies $U_{\ge,\ell}(x,0)+m(\gt,\ge)\geq u(x,\gt)$. If $x \in B_\ge^c$,  $m(\gt,\ge)\geq W_{\ge/2}(x,\gt)\geq u(x,\gt)$, hence  $U_{\ge,\ell}(x,0)+m(\gt,\ge)\geq u(x,\gt)$. So we always have $U_{\ge,\ell}(x,0)+m(\gt,\ge)\geq u(x,\gt)$ for any $x \in {\BBR}^N$. Applying maximum principle yields  to $U_{\ge,\ell}(.,.-\gt)+m(\gt,\ge) \geq u$ in ${\BBR}^N \ti (\gt,\ity)$. Finally, the function $\overline U$ defined by $(\ref{maxs})$ is the maximal solution because $U_{\ge,\ell}(x,t-\gt)\to U_{\ge,\ell}(x,t)$ as $\gt\to 0$ and 
$U_{\ge,\ell}\uparrow U_{\ge,\ity}$ when $\ell\to\infty$ and $U_{\ge,\ity}\downarrow  \overline U$ when $\ge\to0$. \qeda
\mysection{About uniqueness}
We prove first the existence of global radial solutions of $(\ref{A18})$ under the Keller-Osserman condition. \medskip

\noindent{\bf Proof of \rprop{Ex1}.} \smallskip

A solution of $(\ref{A19})$ is locally given by the formula
	\be w(r)=a+\int_0^r s^{1-N}\int_0^s t^{N-1}f(w)dtds \ee
Existence follows from the Picard-Lipschitz fixed point theorem. The function is increasing and defined on a maximal interval $[0,r_a)$. By a result of Vazquez and Veron \cite{VaVe} $r_a=\ity$, thus the solution is global. Uniqueness on $[0,\ity)$ follows always from local uniqueness. The function $r \mt w(r)$ is increasing and 
	\[ w\pr (r)\geq\frac{ah(a)}{N}r, \]
	\[ w(r)\geq a+\frac{a h(a)}{2N}r^2 \]
for all $r>0$. \qeda 
\bprop{Ex2} Assume $(\ref{A16})$ holds. For any $u_0 \in C({\BBR}^N)$ which satisfies
	\be w_a(\abs x)\leq u_0(x) \leq w_b(\abs x) \qq \forall x \in {\BBR}^N	\ee
for some $0<a<b$, there exists a positive function $\overline u \in C(\ovl Q_\ity) \cap C^{2,1}(Q_\ity)$ solution of $(\ref{A1})$ in $Q_\ity$ and satisfying $\overline u(.,0)=u_0$ in $\BBR^N$. Furthermore
	\bel{Ex2.1} w_a(\abs x)\leq u(x,t) \leq w_b(\abs x) \qq \forall (x,t)\in Q_\ity.	\ee
\es
\Proof Clearly $w_a$ and $w_b$ are ordered solutions of $(\ref{A1})$. We denote by $u_n$ the solution of the initial-boundary problem
\be \left\{ \BA{lll}	
\prt_t u_n -\Gd u_n + f(u_n) &= 0 \q  &\text{in } Q_n = B_n \ti (0,\ity) \\\phantom{\prt_t u_n -\Gd u_n ,}
u_n(x,t)&=(w_a(\abs x)+w_b(\abs x))/2 \q &\text{in } \prt B_n \ti (0,\ity) \\\phantom{\prt_t u_n -\Gd u_n ,}
u_n(x,0)&=u_0(x)&\text{in } B_n.
\EA \right. \ee
By the maximum principle, $u_n$ satisfies $(\ref{Ex2.1})$ in $Q_n$. Using locally parabolic equations regularity theory, we derive that the set of functions $\{u_n\}$ is eventually equicontinuous on any compact subset of $\overline Q_\infty$. Using a diagonal sequence, we conclude that there exists a subsequence $\{u_{n_k}\}$ which converges locally uniformly in  $\overline Q_\infty$ to some weak solution $\overline u\in C(\overline Q_\infty)$ which satisfies $\overline u(.,0)=u_0$ in $\BBR^N$. By standard method, $\overline u$ is a strong solution (at least $C^{2,1}(Q_\infty)$).\qeda 
\bprop{Ex3} Assume $(\ref{A16})$ and $(\ref{A10})$ hold. Then for any $u_0 \in C({\BBR}^N)$ which satisfies
	\be w_a(\abs x)\leq u_0(x) \leq w_b(\abs x) \qq \forall x \in {\BBR}^N	\ee
for some $0<a<b$, there exists a positive function $\underline u\in C(\overline{Q_\infty})$ solution of $(\ref{A1})$ in $Q_\ity$ satisfying $\underline u(.,0)=u_0$ in $\BBR^N$ and 
\bel{Ex3.0*} \underline u(x,t) \leq \min\{\gf_\ity(t),w_b(\abs x)\} \qq \forall (x,t) \in Q_\ity. \ee
\es
\Proof For any $R>0$, let $u_R$ be the solution of
\bel{u_R} \left\{ \BA{ll} 
\prt_t u_R - \Gd u_R +f(u_R) = 0 \qq \qq \qq &\text{  in } Q_\ity \\
\phantom{\prt_t u_R - \Gd  +,,}
u_R(x,0)=u_0(x)\gc_{B_R}(x) &\text{ in } {\BBR}^N.
\EA \right. \ee
The solution which is constructed is dominated by the solution of the heat equation with the same initial data. Thus
\bel{Ex3.0} u_R(x,t)\leq (4\gp t)^{-N/2}\myint{B_R}{}e^{-|x-y|^2/4t}u_0(x) dy\qq \forall (x,t) \in Q_\ity.\ee 
and ${\displaystyle \lim_{|x|\to\infty}} u_R(x,t)=0$ uniformly with respect to $t$.
The functions $\gf_\ity$ and $w_b$ are solutions of $(\ref{A1})$ in $Q_\ity$, which dominate $u_R$ at $t=0$. By the maximum principle, 
\bel{Ex3.1} \min\{\gf_\ity(t),w_b(|x|)\} \geq u_R(x,t) \q  \forall (x,t) \in Q_\ity. \ee 
The fact that the mapping $R \mt u_R$ is increasing and $(\ref{Ex3.1})$ imply that there exists $\underline u:={\displaystyle \lim_{R \to \ity}} u_R$ which satisfies $\underline u(.,0)=u_0$ in  ${\BBR}^N$. Letting $R \to \ity$ in $(\ref{Ex3.1})$ yields $(\ref{Ex3.0*})$. \qeda  \medskip

\noindent{\bf Proof of \rth{non-unique}.} \smallskip 

Combining \rprop{Ex2} and \rprop{Ex3}  we see that there exists two solutions $\underline u$ and $\overline u$ with the same initial data $u_0$ which are ordered and different since 
${\displaystyle \lim_{|x|\to\infty}}\overline u(x,t)=\infty$
and
${\displaystyle \lim_{|x|\to\infty}}\underline u(x,t)\leq \gf_\infty(t)<\infty
$ for all $t>0$.
\qeda \medskip

\noindent{\bf Proof of \rth{unique}.} \smallskip

\noindent {\it Step 1: There always holds}
\bel{Uniqueness1.1}  (ah(a)-bh(b))sign(a-b)\geq \abs{a-b}h(\abs{a-b}) \qq \forall a, b>0 \ee
where $h$ is defined in $(\ref{h})$ and 
$$ sign(z)  = \left\{ \BA{ll} 1 \qq &\text{if } z>0, \\
-1 \qq &\text{if } z<0,\\
0 \qq &\text{if } z=0.
\EA \right. $$
In fact, since $h$ is increasing and assuming $a>b$, we get
\[ \BA{l} 
ah(a)-bh(b)=(a-b)h(a)+b(h(a)-h(b))\\
\phantom{ah(a)-bh(b)}
\geq (a-b)h(a) \\
\phantom{ah(a)-bh(b)}
\geq (a-b)h(a-b).
\EA \]
\noindent {\it Step 2: End of the proof}. By Kato's inequality, 
\[ \prt_t \abs {u-\tl u} - \Gd \abs {u-\tl u} \leq [\prt_t (u-\tl u)-\Gd (u- \tl u)]sign(u- \tl u), \]
therefore by step 1,  
\be \prt_t \abs{u-\tl u} -\Gd \abs{u-\tl u}+\abs{u-\tl u}h(\abs{u-\tl u})\leq 0. \ee
Let $\ge>0$. There exists $R_\ge >0$ such that for any $R\geq R_\ge$,
\be 0\leq  \abs {u-\tl u}(x,t) \leq w_\ge (\abs x) \q \forall (x,t) \in B_R^c \ti [0,1]. \ee
Since $w_\ge$ is a positive solution of $(\ref{A1})$ which dominates $ \abs {u-\tl u}$ on $\prt B_R \ti [0,1]$ and at $t=0$, it follows that $ \abs {u-\tl u}\leq w_\ge$ in $B_R \ti [0,1]$. Letting $R \to \ity $ yields to $ \abs {u-\tl u} \leq w_\ge $ in $\BBR^N \ti [0,1]$. Letting $\ge \to 0$ and since ${\displaystyle \lim_{\ge \to 0}}w_\ge(\abs x)=0$ for any $x \in {\BBR}^N$, we derive $ \abs {u-\tl u}=0$, thus $u=\tl u$ in $\BBR^N \ti [0,1] $. Iterating yields that equality holds in $Q_\ity$. \qeda  \medskip

\noindent{\it Remark.} If we replace the condition $(C1)$ by the condition $(C2)$, the conclusion of \rth{unique} remains valid. Indeed, it follows by the convexity of $f$ that
$$ (f(a)-f(b))sign(a-b) \geq f(\abs{a-b}) \forevery a,b>0. $$ 
Then we proceed as in step 2 to get the desired conclusion. \medskip 

\noindent{\bf Proof of \rth{unique-KO}.} \smallskip

\noindent{\it Proof of statement (i).} The solution $\unl u$ which is constructed in \rprop{Ex3} is a minimal solution of $(\ref{A1})$ in $Q_\ity$ with the initial value $u_0$. Indeed, if $u \in C^{2,1}(Q_\ity)$ is a nonnegative solution of $(\ref{A1})$ in $Q_\ity$ which satisfies $u(.,0)=u_0$ in $\BBR^N$ then, by maximum principle, $ u_R \leq u$ in $Q_\ity$ where $u_R$ is the solution of $(\ref{u_R})$. Letting $R \to \ity$ yields $\unl u \leq u$ in $Q_\ity$. Next we construct the maximal solution. We recall that $w_R$ is the solution of $(\ref{w_R})$. Since $f$ is convex, $f'$ is nondecreasing and $w_R$ there holds $ f'(u_R)\leq f'(w_R+u_R)$, thus there holds
$f(w_R)+f(u_R)\leq f(w_R+u_R)$. Consequently $w_R+u_R$ is a supersolution in $B_R \ti (0,\ity)$. If $u \in C(\ovl {Q_\ity})$ is a solution $(\ref{A1})$ in $Q_\ity$ with initial data $u_0$, it is dominated by $w_R+u_R$ on $\prt B_R \ti (0,\ity)$. Thus $u\leq w_R+u_R$.
which dominates any solution $u \in C(\ovl {Q_\ity})$ at 
 of $(\ref{A1})$ in $B_R \ti (0,\ity)$. Since
 $$u_R\leq u\leq w_R+u_R,
 $$
$w_R \to 0$ when $R\to\infty$, by \rprop {Up2}-Step 1, and $u_R\to \underline u$, we derive that $u=\underline u$.
 \medskip

\noindent{\it Step 2: Construction of a minimal solution.} Assume there exists at least one positive solution $u$ of $(\ref{A1})$ satisfying $(\ref{A22'})$ and $f(u)\in L^1_{loc}\overline{Q_\infty})$. equivalently \cite{MV3}
\bel{weak}
\myint{0}{\infty}\myint{\BBR^N}{}\left(-u(\prt_t\eta+\Gd\eta)+f(u)\eta\right)dxdt=\myint{\BBR^N}{}\eta(x,0)d\gm(x)
\ee
for all $\eta\in C_c^{2,1}(\overline{Q_\infty})$. We construct first a minimal solution in the following way: let $n\in\BBN$ and  $R>0$ and let $v=v_{R,n}$ be the solution of 
\bel{X1}\left\{\BA {ll}
\prt v-\Gd v+f(v)=0\qquad&\text {in }B_R \ti (0,\ity)\\[2mm]\phantom{\prt v-\Gd v+f()}
v=0\qquad&\text {in }\prt B_R \ti (0,\ity)\\[2mm]\phantom{-\Gd u+f)}
v(.,0)=u(.,2^{-n})\qquad&\text {in } B_R.
\EA\right.\ee
By the maximum principle, $v_{R,n}(.,t)\leq u(.,t+2^{-n})$. Furthermore, 
$$v_{R,n}(x,2^{-n})\leq u(.,2^{-n+1})=v_{R,n}(x,0),$$ therefore, 
\bel{X2}v_{R,n}(x,t+2^{-n})\leq v_{R,n-1}(x,t)\qquad\text{in }B_R \ti (0,\ity).
\ee
Using the formulation $(\ref{weak})$ with $v_{R,\ge}$, we obtain
\bel{weak1}
\myint{0}{\infty}\myint{\BBR^N}{}\left(-v_{R,n}(\prt_t\eta+\Gd\eta)+f(v_{R,n})\eta\right)dxdt=\myint{\BBR^N}{}\eta(x,0)u(x,2^{-n})dx,
\ee
for any $\eta\in C_c^{2,1}(\overline{Q^{B_R}_\infty})$. The right-hand side of $(\ref{weak1})$ converges to $\int_{\BBR^N}{}\eta(x,0)d\gm(x)$. Concerning the left-hand side, there holds $f(v_{R,n}(x,t))\leq f(u(x,t+2^{-n}))$. Since $f(u)\in L^1_{loc}\overline{Q_\infty})$, $f(v_{R,n})$ is bounded in $L^1_{loc}\overline{Q^{B_R}_\infty})$. By the $L^1$ regularity theory for parabolic equations (see \cite{MV2} and the references therein), the set of functions $\{v_{R,n}\}$ is locally compact in $L^1_{loc}{Q_\infty})$ and there exists a subsequence $\{n_k\}$ and a function $\underline u_R$ such that $v_{R,n_k}\to \underline u_R$, almost everywhere in $Q^{B_R}_\infty$, and  $\underline u_R\leq u$. Noticing that the sets of functions $\{f(u(.,.+2^{-n}))\}$ and  $\{u(.,.+2^{-n})\}$ are uniformly integrable, we obtain that the two sets $\{f(v_{R,n})\}$ and $\{v_{R,n}\}$ are also uniformly integrable in  $B_R\ti(0,T)$. It follows from Vitali's convergence theorem that, up to a subsequence still denoted by $\{n_k\}$, $v_{R,n_k}\to \underline u_R$ and $f(v_{R,n_k})\to  f(\underline u_R)$ in $L^1(B_R\ti(0,T))$. Letting $n=n_k\to\infty$ in $(\ref{weak1})$ we derive 
\bel{weak2}
\myint{0}{\infty}\myint{\BBR^N}{}\left(-\underline u_R(\prt_t\eta+\Gd\eta)+f(\underline u_R)\eta\right)dxdt=\myint{\BBR^N}{}\eta(x,0)d\gm(x).
\ee
This means that $\underline u_R$ satisfies $\underline u_R\leq u$ and
\bel{X2}\left\{\BA {ll}
\prt \underline u_R-\Gd \underline u_R+f(\underline u_R)=0\qquad&\text {in }B_R \ti (0,\ity)\\[2mm]\phantom{\prt \underline u_R-\Gd \underline u_R+f()}
\underline u_R=0\qquad&\text {in }\prt B_R \ti (0,\ity)\\[2mm]\phantom{-\Gd \underline u_R u+f)}
\underline u_R(.,0)=\chi_{_{B_R}}\gm\qquad&\text {in } B_R.
\EA\right.\ee
If $\tilde u$ is any other nonnegative solution of $(\ref{A1})$ in $Q_{\infty}$ with initial data $\gm$, the same construction of $\tilde v_{R,n}$ solution of $(\ref{X1})$ with initial data $\tilde u(.,2^{-n})$ instead of $u(.,2^{-n})$ converges, up to a subsequence to some $\underline {\tilde u}_R$ which satisfies $\underline {\tilde u}_R\leq \tilde u$ and is solution of problem $(\ref{X2})$. We know from \cite{MV1}, \cite{MV2} that this problem admits at most one solution. Therefore $\underline {\tilde u}_R=\underline {u}_R$, which implies that $\underline {u}_R\leq \tilde u$ in $Q^{B_R}_{\infty}$. Furthermore, in the above construction, we have only used the fact that $\tilde u$ is defined in a domain larger than $Q^{B_R}_{\infty}$ and is nonnegative. Consequently, 
the same comparison applies if we compare $\underline {u}_R$ and $\underline {u}_{R'}$ for $R'>R$ and we obtain
$$ \underline {u}_R\leq \underline {u}_{R'}\qquad\text{in }Q^{B_R}_{\infty}.
$$ 
Put $\underline u=\lim_{R\to\infty}\underline {u}_R$. Using the monotone convergence theorem and a test function $\eta\in C_c^{2,1}(\overline{Q_\infty})$ with compact support in $Q^{B_R}_{\infty}$, we obtain 
\bel{weak2}
\myint{0}{\infty}\myint{\BBR^N}{}\left(-\underline u(\prt_t\eta+\Gd\eta)+f(\underline u)\eta\right)dxdt=\myint{\BBR^N}{}\eta(x,0)d\gm(x).
\ee
from $(\ref{weak2})$. Thus $\underline u$ satisfies $(\ref{A22'})$ and $f(\underline u)\in L^1_{loc}\overline{Q_\infty})$. By construction $\underline u$ is smaller than any other nonnegative solution.
 \medskip

\noindent{\it Step 3: Proof of statement (ii).}  As in the proof of statement (i), we see that, for any $n\in\BBN^*$, there holds $u\leq W_R+v_{R,n}$ in $Q^{B_R}_{\infty}$. Consequently $u\leq W_R+\underline u_{R}$ and letting $R\to\infty$, $u\leq \underline u$. Thus $u=\underline u$.
\qeda
\mysection{Initial trace}
If $\Gw$ is an open domain in $\BBR^N$, we denote by $\GTM(\Gw)$ (resp. $\GTM^b(\Gw)$) the set of Radon measures in $\Gw$ (resp. bounded Radon measures), and by $\GTM_+(\Gw)$ (resp. $\GTM_+^b(\Gw)$) its positive cone. For $T>0$, we set $Q_T^\Gw=\Gw \ti (0,T)$. 
\subsection{The regular part of the initial trace}
In this section we only assume that $f$ is a continuous nonnegative function defined on $\BBR_+$ and that $u$ is a $C^{2,1}$ positive solution of $(\ref{A1})$ in $Q_T$.
\blemma{Init1} Assume $G$ is a bounded $C^2$ domain in $\BBR^N$, $Q_T^{\overline G}:=\overline G\ti (0,T]$ and let $u\in C^{2,1}(Q_T^{\overline G})$ be a positive solution of $(\ref {A1})$ in $Q_T^G$ such that $u, f(u)\in L^1(Q_T^G)$. Then $u\in L^\infty(0,T;L^1(G'))$ for any domain $G'\subset\overline{G'}\subset G$ and there exists  a positive Radon measure $\gm_G$ on $G$ such that
\bel{Tr1} \lim_{t \to 0}\myint{G}{} u(x,t)\gz(x)dx=\myint{G}{} \gz(x)d\gm_G(x)\qq\forall \gz\in C_c(G). \ee
\es
\Proof Let $\gf:=\gf_G$ be the first eigenfunction of $-\Gd$ in $W^{1,2}_0(G)$ with corresponding eigenvalue $\gl_G$. We assume $\gf>0$ in $G$. Then
$$\myfrac{d}{dt}\myint{G}{}u\gf dx+\gl_G\myint{G}{}u\gf dx+\myint{G}{}f(u)\gf\,dx+\myint{\prt G}{}u\gf_{\bf n}dS=0
$$ where $\gf_{\bf n}$ denote the outward normal derivative of $\gf$.
Since  $\gf_{\bf n}<0$, the function
$$t\mapsto e^{\gl_Gt}\myint{G}{}u(x,t)\gf(x) dx-\myint{t}{T}\myint{G}{}e^{\gl_Gs}f(u)\gf dx\,ds$$
 is increasing and
$$\myint{G}{}u(x,t)\gf(x) dx\leq e^{\gl_G(T-t)}\myint{G}{}u(x,T)\gf(x) dx+e^{-\gl_Gt}\myint{t}{T}\myint{G}{}e^{\gl_Gs}f(u)\gf dx\,ds
$$
for $0<t\leq T$. Thus $u\in L^\infty(0,T;L^1(G'))$ for any strict domain $G'$ of $G$. If $\gz\in C_c(G)$, there holds
\bel{Tr1'}\myfrac{d}{dt}\left(\myint{G}{}u(x,t)\gz(x)dx-\myint{t}{T}\myint{G}{}\left(f(u)\gz-u \Gd\gz\right) dx\,ds\right)=0.
\ee
Consequently
\bel{Tr2} 
\lim_{t \to 0}\myint{G}{} u(x,t)\gz(x)dx=\myint{G}{} u(x,T)\gz(x)dx+\myint{0}{T}\myint{G}{}\left(f(u)\gz-u \Gd\gz\right) dx\,ds. \ee
This implies that $u(.,t)$ admits a limit in $\CD'(G)$, and this limit is a positive distribution. Therefore there exists a positive Radon measure $\gm_G$ on $G$ satisfies $(\ref{Tr1})$.\qeda\medskip

\noindent {\bf Proof of \rprop{reg}.} \medskip

It is clear that $\CR(u)$ is an open subset. If $G$ is a strict bounded subdomain of $\CR(u)$, i.e. $\overline G\subset \CR(u)$, there exists a finite number of points $z_j$ ($j=1,...,k$) and $r'_j>r_j>0$ such that $u$, $f(u)\in L^1(Q_T^{B_{r'_j}(z_j)})$ and $\overline G\subset {\displaystyle \cup_{j=1}^k}B_{r_j}(z_j)$. Let $\gm_j=\gm_{B_{r_j}(z_j)}$ the measure defined in \rlemma{Init1}. If $\gz\in C_c(G)$ there exists a partition of unity $\{\eta_j\}_{j=1}^k$  relative to the cover $\{B_{r_j}(z_j)\}_{j=1}^k$ such that $\eta_j\in C^\infty_0(G)$,  $ \text{supp}(\eta_j) \sbs {B_{r_j}(z_j)})$ and ${\displaystyle \gz=\sum_{j=1}^k}\eta_j\gz$. Since
$$\lim_{t\to 0}\myint{B_{r_j}(z_j)}{}u(x,t)(\eta_j\gz)(x) dx=\myint{B_{r_j}(z_j)}{}(\eta_j\gz)(x) d\gm_j(x)\qq\forall j=1,...,k,
$$
there exists a positive Radon measure $\gm$ on $\CR(u)$ satisfying $(\ref{A23})$. Notice also that $u\in L^\infty(0,T;L^1(G))$ for any $G\subset\overline G\subset\CR(u)$.\qeda \medskip

The main problem is to analyse the behaviour of $u$ on the singular set $\CS(u)$.

\subsection{The Keller-Osserman condition holds}
If the Keller-Osserman condition holds, the existence of an initial trace of arbitrary positive solutions of $(\ref{A1})$ is based upon a dichotomy in the behaviour of those solutions near $t=0$. 

\blemma{Trace1} Assume $u$ is a positive solution of $(\ref{A1})$ in $Q_T$ and  $z \in \CS(u)$. Suppose that at least one of the following sets of conditions holds.\smallskip

\noindent (i) There exists an open neighborhood $G$ of $z$ such that $u\in L^1(Q_T^G)$. \smallskip 

\noindent (ii) $f$ is nondecreasing and $(\ref{A12})$ holds. \smallskip  

\noindent Then, for every open relative neighborhood $G'$ of $z$,
\bel{Bl-u}
\lim_{t\to 0}\myint{G'}{}u(x,t) dx=\infty.
\ee
\es
\Proof First, we assume that (i) holds and let $\gz\in C^2_c(G)$, $\gz\geq 0$. Since $z \in \CS(u)$, then for every open relative neighborhood $G'$ of $z$, there holds
\bel{Bl-u1}
\myint{0}{T}\myint{G'}{}f(u)dx\,dt=\infty.
\ee
Since there exists
$$\lim_{t\to 0}\myint{t}{T}\myint{G'}{}u\Gd \gz dx\,dt =L \in \BBR,
$$
it follows from $(\ref{Tr2})$ that 
\bel{Bl-u2}\myint{G'}{} u(x,t)\gz(x)dx=\myint{t}{T}\myint{G'}{}f(u)\gz dxds +O(1),
\ee
which implies $(\ref{Bl-u})$.\smallskip

\noindent Next we assume that $(\ref{A12})$ holds and $u\notin L^1(Q_T^G)$ for every relative neighborhood $G$ of $z$. If there exists an open neighborhood $G\subset \Gw$ of $z$ such that $(\ref{Bl-u})$ does not hold, there exists a sequence $\{t_n\}$ decreasing to $0$ and $0\leq M<\infty$ such that
\bel{Bl-u3}
\sup_{t_n}\myint{G}{}u(x,t_n) dx=M.
\ee
Furthermore, we can always replace $G$ by an open ball $B_R(z)\subset G$. Thus $(\ref{Bl-u3})$ holds with $G$ replaced by $B_R(z)$. Let  $w:=w_R$ be the maximal solution of 
\bel{Bl-u4}\left\{\BA {ll}
-\Gd w+f(w)&=0\qq\text {in  }B_R(z)\\[1mm]
{\displaystyle \lim_{|x-z|\to R}}w(x)&=\infty.
\EA\right.\ee
Let $v:=v_n$ be the solution of 
\bel{Bl-u5}\left\{\BA {lll}
 \prt_tv-\Gd v&=0\qq&\text {in  }B_R(z)\ti (t_n,\infty)\\[1mm]\phantom{ \prt_tv-\Gd}
v&=0\qq&\text {in  }\prt B_R(z)\ti (t_n,\infty)\\\phantom{,,, }
\!v(.,t_n)&=u(.,t_n) \qq&\text {in  } B_R(z).
\EA\right.\ee
Since $v_n\geq 0$, $f(w_R+v_n)\geq f(w_R)$, and  $w_R+v_n$ is a supersolution of $(\ref{A1})$ in $B_R(z) \ti (t_n,T)$. It dominates $u$ on $\prt B_R(z)\ti (t_n,T)$ and at $t=t_n$, thus $u\leq w_R+v_n$ in ${B_R(z)}\ti(t_n,T)$.  We can assume that $u(.,t_n)\to\gn$ for some positive and bounded measure $\gn$ on $B_R(z)$. Therefore
\bel{Bl-u6}
u(x,t)\leq v(x,t)+w_R(x)\qq \text {in  }Q_T^{B_R(z)}
\ee
where $v$ is the solution of 
\bel{Bl-u7}\left\{\BA {lll}
 \prt_tv-\Gd v&=0\qq&\text {in  }Q_\infty^{B_R(z)}\\[1mm]\phantom{ \prt_tv-\Gd}
v&=0\qq&\text {in  }\prt B_R(z)\ti (0,\infty)\\\phantom{,,, }
\!v(.,0)&=\gn \qq&\text {in  } \CD'(B_R(z)).
\EA\right.\ee
Since $v\in L^1(Q_T^{B_R(z)})$ and $w_R$ is uniformly bounded in any ball $B_{R'}(z)$ for $0<R'<R$, we conclude that $u\in L^1(Q_T^{B_{R'}(z)})$, which is a contradiction.\qeda

\bdef{Trace} Assume $f$ is nondecreasing and satisfies $(\ref{A12})$. Let $u\in C^{2,1}(Q_T)$ be a positive solution of $(\ref{A1})$ in  $Q_T$. We say that $u$ possesses an initial trace with regular part $\gm\in\GTM_+(\CR(u))$ and singular part $\CS(u)=\BBR^N\setminus\CR(u)$  if\smallskip

\noindent (i) For any $\gz\in C_c(\CR(u))$,
\bel {Bl-u8}
\lim_{t\to 0}\myint{\CR(u)}{}u(x,t)\gz(x) dx=\myint{\CR(u)}{}\gz(x) d\gm(x).
\ee
(ii) For any open set $G\subset\BBR^N$ such that $G\cap\CS(u)\neq\emptyset$
\bel {Bl-u9}
\lim_{t\to 0}\myint{G}{}u(x,t) dx=\infty.
\ee
\es

\noindent {\bf Proof of \rth{tr+KO}} \medskip

The set $\CR(u)$ and the measure $\gm\in\GTM_+(\CR(u))$  are defined by \rdef{Reg} thanks to \rprop{reg}.  Because $(\ref{A12})$ holds, $\CS(u)=\Gw\setminus\CR(u)$ inherits the property (ii) in \rdef{Trace} because of \rlemma {Trace1} (ii).\qeda\medskip

If $\Gw$ is a bounded domain with a $C^2$ boundary and $\gm\in\GTM_+^b(\Gw)$, we denote by $u_\gm$ the solution of
 \bel{meas}\left\{\BA {l}
  \prt_tu-\Gd u+f(u)=0\qq\text{in }Q_\infty^\Gw\\\phantom{  \prt_tu-\Gd +f(u)}
 u=0\qq\text{in }\prt\Gw\ti (0,\infty)\\\phantom{  \prt_tu-\Gd f,,}
 u(.,0)=\gm\qq\text{in }\CD'(\Gw).
 \EA\right.\ee\smallskip
 
 We recall the following stability result proved in \cite[Th 1.1]{MV2}.
 \blemma {stab}Let $\Gw$ be a bounded domain with a $C^2$ boundary. Assume $f$ is nondecreasing and satisfies $(\ref{A8})$. Then for any $\gm \in \GTM^b(\Gw)$ problem $(\ref{meas})$ admits a unique solution $u_\gm$. Moreover, if $\{\gm_n\}\subset\GTM^b(\Gw)$ converges weakly to $\gm\in\GTM^b(\Gw)$ then $u_{\gm_n}\to u_\gm$ locally uniformly in  $\overline \Gw\ti (0,\infty)$ and in $L^1(Q_T^\Gw)$, and $f(u_{\gm_n})\to f(u_{\gm})$ in $L^1(Q_T^\Gw)$, for every $T>0$.
 \es
 
 \noindent\Remark The result remains true if $\Gw$ is unbounded, with a $C^2$ compact (possibly empty) boundary and the $\gm_n$ have their support in a fixed compact set. In such a case $u_{\gm_n}(x,t)\to 0$ when $|x|\to\infty$, uniformly with respect to $n$ and $t$ since
 \bel{EST}
 \abs{u_{\gm_n}(x,t)}\leq \myfrac{1}{(4\gp t)^{N/2}}\myint{\BBR^N}{}e^{-\abs{x-y}^2/4t}d\abs{\gm_n}(y)\quad\forall (x,t)\in Q_\infty.
 \ee
 \medskip

By \rlemma{stab} and the remark hereafter, for every $y \in \Gw$ and $k>0$, there exists a unique solution $v_{y,k,\Gw}:=v$ to $(\ref{meas})$ with $\gm=k\gd_0$. By comparison principle (see \cite[Prop 1.2]{MV2}) $v_{y,k,\Gw}$ is positive, increases as $k$ increases and depends continuously on $y$. Note that if $\Gw=\BBR^N$, $v_{y,k,\BBR^N}(x,t) := v_{y,k}(x,t)=u_k(\abs{x-y},t))$; furthermore, if  $f$  satisfies $(\ref{A12})$, we recall that $\underline U=\lim_{k\to\infty}u_k$ is the minimal solution of $(\ref{A1})$ in $Q_\infty$ with initial trace $(\{0\},0)$.

\bprop{singular-point} Assume $f$ is nondecreasing and satisfies $(\ref{A8})$ and $(\ref{A12})$. Let $u \in C^{2,1}(Q_\ity)$ is a positive solution of $(\ref{A1})$ in $Q_\ity$ with initial trace $(\CS,\gm)$. Then for every $y \in \CS$,
\bel{singular-point1} \underline U_y(x,t):= \underline U(x-y,t) \leq u(x,t) \ee
in $Q_\ity$. 
\es
\Proof By translation we may suppose that $y=0$. Since $0 \in \CS(u)$, for any $\gh>0$ small enough 
$$ \lim_{t \to 0}\myint{B_\gh}{}u(x,t)dx=\ity. $$
For $\ge>0$, denote $M_{\ge,\gh}=\myint{B_\gh}{}u(x,\ge)dx$. For any $m>m_\gh={\displaystyle \inf_{\gs>0}M_{\gs,\gh}}$ there exists $\ge=\ge(m,\gh)$ such that $m=M_{\ge,\gh}$ and ${\displaystyle \lim_{\gh \to 0}\ge(m,\gh)=0}$. Let $v_{\gh}$ be the solution of the problem
$$ \left\{ \BA{lll}
\prt_t v_{\gh} - \Gd v_{\gh} + f(v_{\gh}) = 0 \qq &\text{in } Q_\infty \\ \phantom{\prt_t v_{\gh} - \Gd v_{} +  }
v_{\gh}(x,0)=u(x,\ge)\gc_{B_\gh} \q &\text{in } \BBR^N
\EA \right. $$
where $\gc_{B_\gh}$ is the characteristic function of $B_\gh$. By the maximum principle $v_{\gh}\leq u$ in 
$\BBR^N \ti (\ge,\ity)$. By \rlemma{stab} and the remark after $v_{\gh}$ converges to $v_{0,m}$ when $\gh$ goes to zero. Letting $m$ go to infinity yields $(\ref{singular-point1})$. \qeda \medskip

\bcor{sup} Under the assumption of \rprop{singular-point}, there exists a minimal positive solution $\underline U_{\CS}$ of $(\ref{A1})$ in $Q_\ity$ with initial trace $(\CS,0)$ in the sense that 
\bel{singular-point2} \underline U_\CS(x,t) \leq u(x,t)\qquad\forall (x,t)\in Q_\infty, \ee
for all positive solution $u\in C^{2,1}(Q_\infty)$ of $(\ref{A1})$ with initial trace $(\CS(u),\gm)$.
\es
\Proof If we set $\underline {\tilde U}_\CS=\sup\{U_y:y\in\CS\}$, 
then $\underline {\tilde U}_\CS$ is a subsolution of $(\ref{A1})$. If $u$ is a positive solution of $(\ref{A1})$ with initial trace $(\CS,\gm)$, then $u\geq \underline {\tilde U}_\CS$ by \rprop{singular-point}. Therefore $u$ is larger than the smallest solution of $(\ref{A1})$ in $Q_\ity$ which is above $\underline {\tilde U}_\CS$. We denote this minimal solution by 
$\underline U_\CS$.\qeda
\medskip

If $\CS$ contains some ball $B_R$ we have a more precise result.
\bprop{sup2} Let $u$ be a positive solution of $(\ref{A1})$ in $Q_\infty$ with initial trace $(\CS,\gm)$. We assume that 
$\CS$ has a non-empty interior, and for $R>0$, we denote by $int_R(\CS)$ the set of $y\in \CS$ such that $\overline B_R(y)\subset int_R(\CS)$. Then for any $R'\in (0,R)$ there holds
\bel{singular-point3} 
\lim_{t\to 0}\myfrac{u(x,t)}{\phi_\infty(t)}=1
\ee
uniformly for $x\in \overline B_{R'}(y)$ and $y\in int_R(\CS)$.
\es
\Proof Let $y\in int_R(\CS)$  and $w(x,t)= u(x,t)+W_R(x-y)$. Then $w$ is a supersolution of $(\ref{A1})$ in $Q_{\infty}^{B_R(y)}$ and $\lim_{t\to 0}w(x,t)=\infty$, uniformly with respect to $x\in B_R(y)$, by $(\ref{singular-point1})$. Then, for any $\ge>0$, there exists $t_\ge>0$ such that $w(x,t)\geq \phi_\infty(\ge)$ in $Q_{t_\ge}^{B_R(y)}$. Since $\phi_\infty((t+\ge)$ remains bounded on $\prt B_R(y)\ti (0\infty)$, it follows by the maximum principle that
$$w(x,t)\geq \phi_\infty((t+\ge)\qquad\forall (x,t)\in Q_{\infty}^{B_R(y)}.
$$
Letting $\ge\to 0$ and using the fact that $W_R(x-y)$ remains uniformly bounded when $|x-y|\leq R'$, we derive
\bel{singular-point4} 
u(x,t)\geq \phi_\infty((t)-K_{R'}\qquad\forall (x,t)\in Q_{\infty}^{B_R'(y)}.
\ee
where $K_{R'}=\max\{W_R(x-y):|x-y|\leq R'\}$. Combining this estimate with (\ref{UNIVE}) yields to (\ref{singular-point3}).\qeda
\medskip

The following convergence lemma is obtained by using the arguments of \rlemma{Init1}
\bprop{initial-convergence} Assume $f$ is nondecreasing and satisfies $(\ref{A8})$ and $(\ref{A12})$. Let $\{u_n\}$ be a sequence of positive solutions of $(\ref{A1})$ in $Q_\ity$ with initial trace $(\CS(u_n),\gm_n)$ such that $u_n \to u$ locally uniformly in $Q_\ity$ and let $A$ be an open subset of $\CR(u_n):=\BBR^N\setminus \CS(u_n)$. Then $u$ is a positive solution of $(\ref{A1})$ in $Q_\ity$, with initial trace denoted by $tr_{\BBR^N}(u)=(\CS,\gm)$.  Furthermore, if  $\gm_n(A)$ remains uniformly bounded, then $A\subset \CR:=\BBR^N\setminus \CS$ and $\chi_{_A}\gm_n\to \chi_{_A}\gm$
weakly.
Conversely, if $A\subset \CR(u)$, then $\gm_{n}(K)$ remains bounded independently of $n$, for every compact set $K\subset A$.
\es
\Proof The fact that $u$ is a positive solution of $(\ref{A1})$ in $Q_\ity$ is standard by the weak formulation of the equation. Assume now that $A\cap \CS\neq\emptyset$. Let $z\in A\cap \CS$ and $R>0$ such that 
$\overline B_R(z)\subset A$. By convexity, $u_n$ is bounded from above  in $Q^{B_R(z)}_\infty$ by $v_n+W_R$, where $v_{n,z}$ satisfies
\bel{Q1} \left\{\BA {ll}
\prt_tv-\Gd v+f(v)=0\qquad&\text{in }Q^{B_R(z)}_\infty\\
\phantom{\prt_tv-\Gd+f(v)}v=0\qquad&\text{in }\prt B_R(z)\ti (0,\infty)
\\
\phantom{\prt_t..+f(v)}v(.,0)=\chi_{_{B_R(z)}}\gm_n\qquad&\text{in } B_R(z),
\EA\right.\ee
and $W_R$ is the maximal solution of $(\ref{Bl-u4})$. We can assume that, up to a subsequence, $\chi_{_{B_R(z)}}\gm_{n_k}\to \gm_z\in\frak M^b_+(B_R(z))$ weakly, thus  $v_{n_k,z}\to v_z$ where $v_z$ is the solution of 
\bel{Q2} \left\{\BA {ll}
\prt_tv-\Gd v+f(v)=0\qquad&\text{in }Q^{B_R(z)}_\infty\\
\phantom{\prt_tv-\Gd+f(v)}v=0\qquad&\text{in }\prt B_R(z)\ti (0,\infty)
\\
\phantom{\prt_t..+f(v)}v(.,0)=\gm_z\qquad&\text{in } B_R(z),
\EA\right.\ee
Therefore
\bel{Q3}u\leq v_z+W_R\qquad\text{in }Q^{B_R(z)}_\infty.
\ee
By \rlemma{stab}, it implies that $u\in L^1(Q^{B_{R'}(z)}_T)$ for any $0<R'<R$. Furthermore, if $(\ref {A8})$ is satisfied, then for any positive constant $k$, $s\mapsto s^{N/2}f(s^{-N/2}+k)\in L^1(0,1)$, thus if $v$ is such that $f(v)\in L^1(Q^{B_{R'}(z)}_T)$, there holds $f(v+k)\in L^1(Q^{B_{R'}(z)}_T)$. In particular, since $f(v_z)\in L^1(Q^{B_{R'}(z)}_T$, and if we take $k=\max\{W_R(x):x\in B_{R'}(z)\}$, we derive that $f(u)\in L^1(Q^{B_{R'}(z)}_T)$, and therefore $z\in\CR$, which is a contradiction; thus $A\subset\CR$.  Next, there exist a subsequence $\{n_k\}$ and a bounded positive measure $\tilde \gm$, with support in $A$ such that $\chi_{_A}\gm_{n_k}\to \tilde\gm$ weakly and suppose $\overline B_R(z)\subset A$. Since $u_{n_k}\leq v_{n_k,z}+k$ and $f(u_{n_k})\leq f(v_{n_k,z}+k)$ in $Q^{B_{R'}(z)}_T$
and $v_{n_k,z}+k$ and $f(v_{n_k,z}+k)$ are uniformly integrable in $Q^{B_{R'}(z)}_T$, it follows that $u_{n_k}$ and 
$f(u_{n_k,z})$ inherit this property. Therefore, if $\gz\in C^2_c(B_R(z))$ we can assume that it vanishes outside $B_{R'}(z)$. Because
\bel{i-c2} \BA {l}
\myint{B_R(z)}{}\gz(x) d\gm_{n_k}(x)= \myint{B_R(z)}{}u_{n_k}(x,t)\gz(x) dx+\myint{0}{t}\myint{B_R(z)}{}
\left(-u_{n_k}\Gd\gz+f(u_{n_k})\gz\right)dxds,\EA\ee
we derive from Vitali's convergence theorem
\bel{i-c3} \BA {l}
\myint{B_R(z)}{}\gz(x) d\tilde \gm(x)= \myint{B_R(z)}{}u(x,t)\gz(x) dx+\myint{0}{t}\myint{B_R(z)}{}
\left(-u\Gd\gz+f(u)\gz\right)dxds.\EA\ee
This implies that $\chi_{_{B_R(z)}}\tilde \gm=\chi_{_{B_R(z)}}\gm$ and, by a partition of unity, that $\tilde \gm=\chi_{_{A}}\gm$.\smallskip

\noindent Assume now that $K\subset\CR$ is compact. If $\gm_n(K)$ is unbounded and up to a subsequence still denoted by $\{n\}$, there exists a point 
$y\in K$ such that for any neighborhood $\CO$ of $y$, $\CO\subset A$, 
$\gm_{n}(\CO)\to\infty$ as $n\to\infty$. We can take $\CO=B_r(y)$ and put 
$M_{n,r}=\gm_{n}(B_r(y))$. If $m\in\BBN^*$, there exists an integer $n=n(m,r)$ such that $m\leq M_{n,r}$, and $\lim_{r\to 0}n(m,r)=\infty$. Let $r_0>r$ such that $B_{r_0}(y)\subset A$, and $w_r$ be the solution of 
\bel{i-c3} \left\{\BA {ll}
\prt_tw-\Gd w+f(w)=0\qquad&\text{in }Q_\infty^{B_{r_0}(y)}\\\phantom{\prt_t-\Gd w+f(w)}
w=0\qquad&\text{in }\prt B_\infty^{B_{r_0}(y)}\\\phantom{\prt_tw-\Gd w+}
\!w(.,0)=\chi_{_{B_r(y)}}\gm_n\qquad&\text{in }B_{r_0}(y).
\EA\right.\ee
By the comparison principle, $w_r\leq u_n$ in  $Q_\infty^{B_{r_0}(y)}$. Since $\chi_{_{B_r(y)}}\gm_n\to m\gd_y$ as $r\to 0$ and $n\to\infty$, we derive $u_{y,m, B_{r_0}(y)}\leq u$ from \rlemma{stab} and the remark hereafter. Since $m$ is arbitrary, $u_{y,\infty, B_{r_0}(y)}\leq u$. This implies that $y\in \CS$, a contradiction.
 \qeda \medskip

If $A$ is an open subset of $\Gw$ and $\gn \in \GTM^+(A)$, we define an extension $\unl \gn$ of $\gn$ to $\Gw$ by
\bel{extension} \unl \gn(E)=\inf_{E \subseteq O}\gn(O \cap A) \ee
for every Borel set $E \sbs \Gw$ where the infimum is taken over the open subsets $O$; $\unl \gn$ is an outer regular Borel measure on $\Gw$ and $\gn=\underline\gn_{\mid A}$. \medskip

The following result which shows the existence of a minimal solution of $(\ref{A1})$ with a given initial trace in $\GTM_+(A)$ for any open subset $A$ in $\BBR^N$ is a straightforward adaptation of \cite[Lemma 3.3]{MV1}.

\bprop{min-restriction} Assume $f$ is nondecreasing and satisfies $(\ref{A8})$, $(\ref{A12})$ and $(C2)$. \medskip

\noindent (i) Let $A$ be an open subset of $\BBR^N$ and let $\gn \in\GTM_+(A)$ with associated extension $\underline\gn$. Then there exists a positive solution of $(\ref{A1})$ in $Q_\ity$ denoted by $\unl u_\gn$ satisfying
$Tr_{\BBR^N}(\unl u_\gn)=\unl \gn $ and such that $\unl u_\gn \leq v$ for every positive solution $v$ of $(\ref{A1})$ in $Q_\infty$ such that $tr_{\BBR^N}(v)=(\CS,\gm)$ and $\chi_{_A}\gm \geq \gn$. \medskip

\noindent (ii) Let $\Gw\subset \BBR^N$ be a bounded domain with a $C^2$ boundary and $u_n$ be the solution of problem 
\bel{U-inf1} \left\{ \BA{lll}
\prt_t u_n -\Gd u_n + f(u_n) &= 0 \qq &\text{in } Q_T^\Gw \\ \phantom{--------}
u_n &= n &\text{on } \prt \Gw \ti (0,\ity) \\ \phantom{------}
u_n(.,0) &= n &\text{in } \Gw.
\EA \right. \ee
Denote $U_{\ity,\Gw}:={\displaystyle \lim_{n \to \ity}u_n}$. Then $U_{\ity,\Gw}$ is the maximal solution of $(\ref{A1})$ in $Q_\ity^\Gw$ in the sense that the following relation holds in $Q_T^\Gw$ for every positive solution $v$ of $(\ref{A1})$
\bel{U-inf2} U_{\ity,\Gw} \geq v. \ee
\es

Taking $A=\CR:=\BBR^N\setminus\CS$, we obtain the existence of a minimal positive solution of $(\ref{A1})$ with a given positive Radon measure $\gm\in\mathfrak M_+(\CR)$ as the regular part of the initial trace.

\bcor{min-restriction1} Let $\CS$ be a closed subset of $\BBR^N$, $\CR=\BBR^N\setminus\CS$ and $\gm\in\mathfrak M_+(\CR)$. Then there exists a positive solution $\underline u_\gm$ of $(\ref{A1})$ such that
 $Tr_{\BBR^N}(\underline u_\gm)=\unl \gm $  and $\unl u_\gm \leq v$ for every positive solution $v$ of $(\ref{A1})$ in $Q_\infty$ such that $tr_{\BBR^N}(v)=(\CS,\gm)$.
\es
\medskip

As a counterpart of \rth{tr+KO} we have the following existence theorem. \medskip

\noindent{\bf Proof of \rth{KO-existence}} \smallskip 

\noindent {\it Step 1: Construction of a minimal solution.} 
Let $\underline u_{\CS}$ and $\underline u_{\gm}$ the minimal solution constructed in \rcor{sup} and \rcor{min-restriction1}. Then $\underline{\check u}_{\CS,\gm}:=\sup\{\underline u_{\CS},\underline u_{\gm}\}$ is a subsolution of $(\ref{A1})$ in $Q_\infty$ while $\underline{\hat u}_{\CS,\gm}:=\underline u_{\CS}+\underline u_{\gm}$ is a supersolution. Furthermore $\underline{\check u}_{\CS,\gm}\leq \underline{\hat u}_{\CS,\gm}$. Therefore the set of solutions $u$  in $Q_\infty$ such that   $\underline{\tilde u}_{\CS,\gm}\leq u\leq \underline{\hat u}_{\CS,\gm}$ is not empty and we denote by $\underline{u}_{\CS,\gm}$ the smallest solution larger than  $\underline{\check u}_{\CS,\gm}$; it is a solution with initial trace $(\CS,\gm)$. If $u$ is any other positive solution with the same initial trace, it is larger than $\underline u_{\CS}$ and $\underline u_{\gm}$ by \rcor{sup} and \rcor{min-restriction1}. Therefore it is larger than $\underline{\check u}_{\CS,\gm}$ and consequently larger than $\underline{u}_{\CS,\gm}$.

\noindent{\it Step 2: Construction of the maximal solution.} 
The proof is somewhat similar to the one on \cite[Th 3-4]{MV1}, but we give it for the sake of completeness. We denote, for $\gd>0$, 
$$\CS^\gd:=\{x \in \BBR^N: dist(x,\CS) \leq \gd\}\;\text{and }\CR^\gd:=\BBR^N \sms \CS^\gd.$$
and let $\gm_\gd$ be the measure given by
$$ \gm_\gd(E)=\gm(\CR_\gd \cap E) \forevery E \sbs \BBR^N, E \text{ Borel}. $$ 
We denote by $u_{\CS^\gd}$ a solution of (\ref{A1}) in $Q_\infty$ with initial trace $(\CS^\gd,0)$: a solution is easily constructed as the limit when $R,k\to\infty$ of the solution $v=v_{k,R}$ of
\bel{Y1}\left\{\BA{l}
\prt_t v-\Gd v+f(v)=0\qquad\text{in }Q_\infty\\\phantom{\prt_t v-\Gd v.}
v(.,0)=k\chi_{_{(\overline B_R\cap\CS^\gd)\cup(\overline B_R\cap \overline B^c_{R-\gd})}}
\EA\right.\ee
By \rprop{sup2}, there holds,  for any $0<\gd'<\gd$ and $\ge>0$, 
\bel{Y2}
\lim_{t\to 0}\myfrac{u_{\CS^\gd}(x,t)}{\phi_\infty(t)}=1\quad\text {uniformly on }\CS_{\gd'}
\ee
 
Let $u_{\gm_{\gd}}$ be the solution of of (\ref{A1}) in $Q_\infty$ with initial trace $(\emptyset,\gm_{\gd})$. This solution is constructed by approximation, as the limit, when $R\to\infty$, of the solution $u=u_{\chi_{_{ B_R}}\gm_{\gd}}$ of 
\bel{Y3}\left\{\BA{ll}
\prt_t u-\Gd u+f(u)=0\qquad&\text{in }Q_\infty\\\phantom{\prt_t v-\Gd v...}
\!u(.,0)=\chi_{_{ B_R}}\gm_{\gd}\qquad&\text{in }\BBR^N.
\EA\right.\ee
For $\gt>0$, let $u_{\gd,\gt}$ be the solution of $(\ref{A1})$ in $Q_\infty$ with initial data $m_{\gd,\gt}$ defined by
$$m_{\gd,\gt}(x)=\left\{\BA{ll}
\phi_\infty(\gt)\qquad&\text{if }x\in \CS_\gd\\
u_{\gm_{\gd}}(x,\gt)\qquad&\text{if }x\in \CR_\gd.
\EA\right.$$
Then $u(.,\gt)\leq m_{\gd,\gt}$ in $\CS_\gd$ and $u(.,\gt)\geq m_{\gd,\gt}$  in $\CR_\gd$ by \rprop{min-restriction}. Therefore
$$\lim_{\gt\to 0} (u(.,\gt)-m_{\gd,\gt}(.))_+=0
$$
in the weak sense of measures. Furthermore, this solution does not depend on $u$, by only on $\CS_\gd$ and $\gm_\gd$. The set of functions $\{u_{\gd,\gt}\}_{\gt>0}$ is locally uniformly bounded in $Q_\infty$. By the 
regularity theory for parabolic equations, there exists a subsequence $\{\gt_k\}$ and a positive solution $u^*_\gd$ of $(\ref{A1})$
in $Q_\infty$ such that $u_{\gd,\gt_k}\to u^*_\gd$ locally uniformly in $Q_\infty$. By \rprop{sup2} and \rprop{min-restriction}, 
$tr_{\BBR^N}(u^*_\gd)= (\CS^\gd,\gm_\gd)$.
Let $\gw_{\gd,\gt}$ be the solution of $(\ref{A1})$ in $Q_\infty$ with initial data $(u(.,\gt)-m_{\gd,\gt}(.))_+$ (it is constructed in the same way as $\underline u_\gm$ in \rprop{min-restriction} -(i)). By \rth{unique-KO}-(ii), $\lim_{\gt\to 0}\gw_{\gd,\gt}=0$, locally uniformly. Since $u\leq u_{\gd,\gt}+\gw_{\gd,\gt}$ in $(\gt,\infty)\ti\BBR^N$, we obtain $u\leq u^*_\gd$. If $0<\gd'<\gd$, we can compare similarly $u_{\gd,\gt}$ with the solution  $u_{\gd',\gt}$ of $(\ref{A1})$ with initial data
$$m_{\gd',\gt}(x)=\left\{\BA{ll}
\phi_\infty(\gt)\qquad&\text{if }x\in \CS_\gd'\\
u_{\gm_{\gd'}}(x,\gt)\qquad&\text{if }x\in \CR_\gd'.
\EA\right.$$
If $u^*_{\gd'}$ is the limit of any sequence  $\{u_{\gd',\gt_{k'}}\}$, it satisfies $0<u^*_{\gd'}\leq u^*_\gd$ and has initial trace $(\CS^{\gd'},\gm_{\gd'})$. If take in particular  $\gd=\gd_n=2^{-n}$, we construct a decreasing sequence of positive solutions $\{u^*_{2^{-n}}\}$ of $(\ref{A1})$ in $Q_\infty$, with $tr_{\BBR^N}(u^*_{2^{-n}})= (\CS^{2^{-n}},\gm_{2^{-n}})$, satisfying 
$$u\leq u^*_{2^{-n}}\qquad\text{in }Q_\infty.
$$
Clearly the limit $\overline u_{\CS,\gm}$ of the sequence $\{u^*_{2^{-n}}\}$ is a positive solution of $(\ref{A1})$ in $Q_\infty$
 with initial trace $(\CS,\gm)$ and is independent of $u$. It is the maximal solution of the equation with this initial trace.
\qeda
\medskip

\noindent \Remark When $f(r)=|r|^{q-1}r$ with $1<q<1+2/N$, precise expansion of $u_{\infty\gd}(x,t)$, when $t\to 0$ allows to prove uniqueness. Even when $f(r)=r\ln^\ga(r+1) $ with $\ga>2$, uniqueness is not known. The first step would be to prove that uniqueness holds if $tr_\Gw(u)=(\{a\},0)$ for some $a\in\Gw$. However, if 
$\CS=\emptyset$, uniqueness holds from \rth{unique-KO}-(ii).

\subsection{The Keller-Osserman condition does not hold}
 In this section we assume that $(\ref{A12})$ does not hold but  $(\ref{A8})$ is satisfied.    
 \blemma {Jfinite}Assume $(\ref{A10})$, $(\ref{A16})$ are satisfied and $\displaystyle\lim u_{k\to\infty}u_k=\phi_\infty$. If $u$ is a positive solution of $(\ref{A1})$ in $Q_{\infty}$ which satisfies
  \bel{J-1}
  \limsup_{t\to 0}\myint{G}{}u(x,t)dx=\infty,
  \ee
  for some bounded open subset $G\subset\BBR^N$, then $u(x,t)\geq \gf_\infty(t)$. This holds in particular if $(C1)$ and $(C3)$  are satisfied.
 \es
 \Proof By assumpion, there exists a sequence $\{t_n\}$ decreasing to $0$ such that 
   \bel{J-2}
  \lim_{n\to \infty}\myint{G}{}u(x,t_n)dx=\infty.
  \ee
 If $(\ref{J-1})$, we can construct a decreasing sequence of open subsets $G_k\subset G$ such that $\overline {G_k}\subset G_{k-1}$, diam$(G_k)=\ge_k\to 0$ when $k\to\infty$, and
    \bel{J-3}
  \lim_{n\to \infty}\myint{G_k}{}u(x,t_n)dx=\infty\qq\forall k\in\BBN.
  \ee
Furthermore there exists a unique $a\in\cap_k G_k$. We set 
$$\myint{G_k}{}u(x,t_n)dx=M_{n,k}.
$$
Since ${\dsps \lim_{n\to\infty}M_{n,k}=\infty}$, we claim that for any $m>0$ and any $k$, there exists $n=n(k)\in\BBN$ such that 
    \bel{J-4}\myint{G_k}{}u(x,t_{n(k)})dx\geq m.
  \ee
By induction, we define $n(1)$ as the smallest integer $n$ such that $M_{n,1}\geq m$. This is always possible. Then we define  $n(2)$ as the smallest integer larger than $n(1)$  such that $M_{n,2}\geq m$. By induction, $n(k)$ is the smallest integer $n$ larger than $n(k-1)$ such that $M_{n,k}\geq m$. Next, for any $k$, there exists $\ell=\ell (k)$ such that
    \bel{J-5}\myint{G_k}{}\inf\{u(x,t_{n(k)});\ell\}dx= m
  \ee
  and we set 
  $$V_{k}(x)=\inf\{u(x,t_{n(k)});\ell\}\chi_{_{G_k}}(x).
  $$ 
  Let $v_k=v$ be the unique bounded solution of 
      \bel{J-6}\left\{\BA {ll}
      \prt_tv-\Gd v+f(v)=0\qq&\text{in }Q_{\infty}\\\phantom{,, v+f(v)}
      v(.,0)=V_{k}\qq&\text{in }\BBR^N.
 \EA\right. \ee
 Since $v(x,0)\leq u(x,t_{n(k)})$, we derive      
  \bel{J-7} 
  u(x,t+t_{n(k)}) \geq v_k(x,t)\qq\forall (x,t)\in Q_\infty.
  \ee
  When $k\to\infty$, $V_{k}\rightarrow m\gd_a$, thus $v_k\to u_{m\gd_a}$ by \rlemma {stab}. Therefore $u\geq u_{m\gd_a}$. Since $m$ is arbitrary and $ u_{m\gd_a}\to\phi_\infty$ when $m\to\infty$ by \rth{TTA}, it follows that 
  $u\geq \phi_\infty$.\qeda
  
   \blemma {Jinfinite}Assume $(\ref{A15})$ and $\displaystyle\lim u_{k\to\infty}u_k=\infty$ hold. There exists no positive solution $u$ of $(\ref{A1})$ in $Q_{\infty}$ which satisfies $(\ref {J-1})$
for some bounded open subset $G\subset\BBR^N$. This holds in particular if $(C1)$ and $(C3)$ are satisfied
 \es
 \Proof If we assume that such a $u$ exists, we proceed as in the proof of the previous lemma. Since \rlemma{stab} holds, we derive that $u\geq u_{m\gd_a}$ for any $m$. Since ${\dsps \lim_{m\to\infty}u_{m\gd_a}(x,t)=\infty}$ for all $(x,t)\in Q_\infty$, we are led to a contradiction.\qeda\medskip
 
 Thanks to these results, we can characterize the initial trace of positive solutions of $(\ref{A1})$ when the Keller-Osserman condition does not hold. \medskip
 
\noindent {\bf Proof of \rth{tr+J-fin}} \medskip

If there exists some open subset $G$ of $\BBR^N$ with the property $(\ref{J-1})$, then $u\geq \phi_\infty$ and the initial trace of $u$ is the Borel measure $\gn_\infty$. Next we assume that for any bounded open subset $G$ of $\BBR^N$ there holds
  \bel{J-8}
  \limsup_{t\to 0}\myint{G}{}u(x,t)dx<\infty.
  \ee 
If $\CS(u)\neq\emptyset$, there exist $z\in \BBR^N$ and an bounded open neighborhood $G$ of $z$ such that 
$$\myint{0}{T}\myint{G}{} f(u) dx dxt=\infty.
$$
By $(\ref{J-8})$, $u\in L^\infty(0,T;L^1(G))\subset L^1(Q_T^G)$. Then, by \rlemma {Trace1}, $(\ref{Bl-u})$ holds, which contradict $(\ref{J-8})$. Thus $\CS(u)=\emptyset$ and $\CR(u)=\BBR^N$. It follows from \rprop{reg} that there exists a positive Radon measure $\gm$ such that 
  \bel{J-9}
  \lim_{t\to 0}\myint{\BBR^N}{}u(x,t)\gz(x)dx=\myint{\BBR^N}{}\gz(x)d\gm(x)\qq\forall\gz\in C_c(\BBR^N).
  \ee
  \qeda\medskip
  
  Because of the lack of uniqueness from \rth {non-unique} it is difficult to give a complete characterization of admissible initial data for solutions of $(\ref{A1})$ under the assumptions of \rth {tr+J-fin}. However, we have the result as in \rprop{exit-2}. \medskip

\noindent{\bf Proof of \rprop{exit-2}} \smallskip

We first notice that $\max\{\phi_\infty(t);w_b(|x|)\}$ is a subsolution of $(\ref{A1})$ which is dominated by the supersolution $\phi_\infty(t)+w_b(|x|)$. The construction is standard: for $\gt>0$ we set 
  $$\psi(x,\gt)=\myfrac{1}{2}\left(\max\{\phi_\infty(t);w_b(|x|)\}+\phi_\infty(t)+w_b(|x|)\right).
  $$
  There exists a function $u=u_\gt\in C(\overline{Q_{\infty}})$ solution of $(\ref{A1})$ in $Q_{\infty}$ satisfying 
  $u_\gt(.,0)=\psi(.,\gt)$. Furthermore
      \bel{J-11}
\max\{\phi_\infty(t+\gt);w_b(|x|)\}\leq  u_\gt(x,t)\leq \phi_\infty(t+\gt)+w_b(|x|)\qquad \forall (x,t)\in Q_\infty.
  \ee
  By the parabolic equation regularity theory, the set $\{u_\gt\}_{\gt>0}$ is locally equicontinuous in $Q_\infty$. Thus there exist a subsequence $\{\gt_n\}$ and $u\in C(Q_\infty)$ such that $u_{\gt_n}\to u$ on any compact subset of $Q_\infty$. Clearly $u$ is a weak, thus a strong solutions of $(\ref{A1})$ and it satisfies ($\ref{J-10}$). Since any solution $u$ with initial trace $\gn_\infty$ dominates $\phi_\infty$ by \rlemma {Jfinite}, it follows that $\phi_\infty$ is the minimal one.\qeda \medskip
  
\noindent {\bf Proof of \rth{tr+J-infin}} \medskip

As in the proof of  \rth{tr+J-fin}  and because of \rlemma{Jinfinite}, $\CS(u)=\emptyset$. Therefore $\CR(u)=\BBR^N$ and the proof follows from \rprop{reg}.\qeda\medskip
 
 \noindent\Remark Under the assumptions of \rth{tr+J-fin}, it is clear, from the proof of  \rprop{Ex2}, that for any $0<a<b$ and any initial data $u_0\in C(\BBR^N)$ satisfying 
 $$w_a(x)\leq u_0(x)\leq w_b(x)\qq\forall x\in\BBR^N
 $$
 there exists a solution $u\in C(\overline {Q_\infty})$ of  ($\ref{A1}$) in $Q_\infty$ satisfying $u(.,0)=u_0$ and 
 $$w_a(x)\leq u(x,t)\leq w_b(x)\qq\forall (x,t)\in Q_\infty.
 $$
 We conjecture that for any positive measure $\gm$ on $\BBR^N$  which satisfies, for some $b>0$, 
       \bel{J-12}
\myint{B_R}{}d\gm(x)\leq \myint{B_R}{}w_b(x) dx\qq\forall R>0
  \ee
  there exists a positive solution $u$ of $(\ref{A1})$ in $Q_\infty$ with initial trace $\gm$. Another interesting open problem is to see if there exist local solutions in $Q_T$ with an initial trace $\gm$ satisfying
         \bel{J-13}
\lim_{R\to\infty}\myfrac{\myint{B_R}{}d\gm(x)}{ \myint{B_R}{}w_b(x) dx}=\infty\qq\forall b>0.
  \ee

\end{document}